\def\be{\begin{eqnarray}}
\def\en{\end{eqnarray}}
\def\UB{\begin{equation}}
\def\UE#1{\label{#1}\end{equation}}
\def\BE{\begin{equation}}
\def\EE#1{\label{#1}\end{equation}}
\def\C{{\mathbb C}}
\def\R{{\mathbb R}}
\def\T{{\mathbb T}}
\def\subs#1#2{\mbox{\small${#1\atop#2}$}}
\def\nL{\nabla^{\rm L}}
\def\rf#1{(\ref{#1})}
\def\i{{\rm i}}
\def\e{{\rm e}}
\def\d{{\rm d}}
\def\bv{{\bm v}}
\def\bU{{\bm U}}
\def\bx{{\bm x}}
\def\by{{\bm y}}
\def\bxi{{\bm\xi}}
\def\ba{{\bm a}}
\def\bom{{\bm\omega}}
\def\oz{{\bm\omega}^{\rm(init)}}
\def\vz{{\bm v}^{\rm(init)}}
\def\phig{\varphi_{\mathrm{g}}}
\def\pi{\uppi}
\documentclass{jfm}
\usepackage{amsmath,amsfonts,amssymb,upmath,graphicx,bm,natbib,url,upgreek}

\title{Time-analyticity of Lagrangian particle trajectories in ideal fluid flow}
\author[V.~Zheligovsky and U.~Frisch]{
V\ls L\ls A\ls D\ls I\ls S\ls L\ls A\ls V\ns Z\ls H\ls E\ls L\ls I\ls G\ls O\ls V\ls S\ls K\ls Y$^{1,2}$\ns \and\ns
U\ls R\ls I\ls E\ls L\ns F\ls R\ls I\ls S\ls C\ls H$^2$}
\affiliation{
$^1$Institute of Earthquake Prediction Theory and Mathematical Geophysics
of the Russian Academy of Sciences, 84/32 Profsoyuznaya St., 117997 Moscow,
Russian Federation
\\[\affilskip]
$^2$UNS,~CNRS,~Lab.~Lagrange,~OCA,~CS~34229,~06304~Nice~Cedex~4,~France}
\begin{document}
\maketitle

\begin{abstract}
It is known that the Eulerian and Lagrangian structures of fluid flow
can be drastically different; for example, ideal fluid flow can have
a trivial (static) Eulerian structure, while displaying chaotic
streamlines. Here we show that ideal flow with limited spatial
smoothness (an initial vorticity that is just a little better than
continuous), nevertheless has time-analytic Lagrangian trajectories
before the initial limited smoothness is lost. For proving such
results we use a little-known Lagrangian formulation of ideal fluid
flow derived by Cauchy in 1815 in a manuscript submitted for a prize
of the French Academy. This formulation leads to simple recurrence
relations among the time-Taylor coefficients of the Lagrangian map
from initial to current fluid particle positions; the coefficients can
then be bounded using elementary methods. We first consider
various classes of incompressible fluid flow, governed by the Euler
equations, and then turn to highly compressible flow governed by the Euler--Poisson
equations, a case of cosmological relevance. The recurrence relation
associated to the Lagrangian formulation of these incompressible and the
compressible problems are so closely related that the proofs of
time-analyticity are basically identical.
\end{abstract}

\section{Introduction}\label{s:intro}

\citet{Tay} in a paper on the spectrum of turbulence observed that ``If
the velocity of the air stream which carries the eddies is very much greater
than the turbulent velocity, one may assume that the sequence of changes
in $u$ [the velocity] at the fixed point are simply due to the passage
of an unchanging pattern of turbulent motion over the point, ...'' This
{\it Taylor hypothesis} is used by experimentalists working with wind
tunnels to obtain approximate information about the spatial structure
of turbulent flow from the temporal variation at the fixed location of a probe.

Taylor's argument has also implications for the spatio-temporal
smoothness in Eulerian coordinates: near a space-time location of
non-vanishing velocity $\bU$ the temporal regularity cannot be
better than the spatial regularity along the direction of $\bU$, as
can be seen, at least heuristically, by performing a Galilean
transformation with velocity $\bU$. For example, if the spatial
variation of the velocity is roughly as the cubic root of the distance
(H\"older continuity of exponent 1/3), then so can be expected for the
temporal variation. Recently obtained weak (distributional) solutions
of the three-dimensional (3D) incompressible Euler equations give precisely
the same spatial and temporal regularity in Eulerian coordinates \citep{BuLeSz},
and the bounds for the material derivative $\partial_t+(\bv\cdot\nabla)$
are better than the bounds for the spatial and time derivatives \citep{Is}.
In Lagrangian coordinates, where we are following fluid particles, the situation
is very different. As we shall see, with sufficiently regular initial
conditions, the Lagrangian temporal
smoothness vastly exceeds the spatial smoothness.

It is well known that some major mathematical questions remain open for 3D
incompressible flow, be it viscous (Navier--Stokes) or ideal/inviscid (Euler).
However, a number of theorems state that, when the initial
data are sufficiently regular, solutions with this regularity
exist for some time, and are unique. All such theorems have been proved
under the assumption that the initial vorticity is slightly more regular
than merely continuous --- for example, that the initial vorticity is
H\"older-continuous, i.e., it is bounded and its spatial increments
between any two points are bounded by a positive fractional power of the
distance between these points. Well-known instances of such results are
the all-time regularity
of 2D flow \citep{Wo33,Ho33} and the finite-time regularity of
3D flow \citep{Li27}. Their proofs make use of the boundedness of the
velocity gradient, which is not guaranteed, even initially, when the
vorticity is just continuous (such matters are discussed further in
Section~\ref{ss:holder}).

Within such and similar frameworks, we shall show in this paper that,
while the assumed initial regularity persists, the Lagrangian fluid
particle trajectories are analytic functions of the time.

As usual, there were precursors --- foremost, and nearly two hundred years ago,
\citet{Cau27}. His Lagrangian formulation of the 3D incompressible
Euler equations, in terms of Lagrangian invariants, provides us with a simple
framework within which we can develop an elementary theory of analytic fluid
particle trajectories. Cauchy's work on invariants has been almost completely
forgotten in the twentieth century and little cited in the nineteenth.
\citet{Han} and \citet{St48,St83} are notorious exceptions; a detailed review
of the history of Cauchy's half-forgotten result can be found in \citet{FrVi}.

\citet{Sto}, in the last chapter of his book on water waves, where he
discussed the Lagrangian theory of free-surface flow, presented
some of the work of his Courant Institute student Frederick Pohle.
This work remains unpublished except in the PhD thesis \citep{Phl}. Pohle
rediscovered the Lagrangian invariant formulation of Cauchy in the 2D
case. He was not aware of Cauchy's 3D work, but pointed out
briefly that his theory can be extended to three
dimensions. Pohle noticed that this Lagrangian formulation leads to
recurrence relations among the formal time-Taylor
coefficients. He wrote them all, but used only a few, since his goal
was to apply the theory to the breakup of dams. Pohle did not
prove the convergence of this formal Taylor series, i.e.~analyticity,
but \citet{Sto} conjectured that analyticity can be proved using
techniques previously used by Lichtenstein. That is precisely what we
are doing here more than half a century later.

More recently, \citet{Se92,Se95a,Se95b}
argued that weak spatial regularity can nevertheless be consistent
with time-analyticity of Lagrangian trajectories. He proposed studying
this using a Biot-Savart-type integral formulation of the Lagrangian
equations and then making use of the mathematical theory of analytic
functionals (a kind of infinite-dimensional generalisation of
functions of many complex variables). \citet{Sh12} tackled the problem,
avoiding the Biot-Savart integrals. These papers by Serfati and Shnirelman are
very interesting for those with a taste for advanced mathematics.
Temporal analyticity of ``particle trajectories'' was also proved in the recent
paper \citep{CVW} for equations of hydrodynamic type with the use of rather
involved combinatorial relations (called {\it ibid.} ``magic identities'').
Such abstract proofs do not provide us with concrete ways of constructing
the time-Taylor series for the motion of fluid particles.

\citet{FriZh} gave an elementary
proof of analyticity that uses Cauchy's Lagrangian formulation. It was
then realised that our method can be generalised, for
example, to flow that is initially analytic in the space
variables. Furthermore, as conjectured by \citet{Se95a,Se95b},
analyticity of Lagrangian fluid particle trajectories holds also for
compressible flow. The present paper addresses all these
issues. Intended primarily for fluid dynamicists, it is self-contained
and does not rely on advanced mathematical techniques.

The outline of the paper is as follows. Section~\ref{s:ideal} is entirely devoted
to ideal incompressible fluid flow. In Section~\ref{ss:cauchyinv} we present the
Lagrangian formulation of the ideal fluid equations due to Cauchy. In
Section~\ref{ss:taylor} we show how it yields recurrence
relations for time-Taylor coefficients of the Lagrangian map. We then
prove the analyticity of fluid-particle trajectories for the simplest
case: an initial vorticity field whose Fourier series is absolutely
convergent (Section~\ref{ss:L1}). The same technique is then adapted to
initial vorticities that are analytic in space (Section~\ref{ss:analytic}) or
H\"older-continuous (Section~\ref{ss:holder}). In Section~\ref{ss:bounds}
we bracket the radius of convergence of these time-Taylor series more
precisely: by considering an example of a Beltrami flow with integrable
Lagrangian trajectories and by refining our previously obtained bounds.
Section~\ref{s:cosmo} is devoted to an instance of compressible ideal
fluid flow of cosmological relevance.
Concluding remarks are made in Section~\ref{s:conclusion}.

\section{Ideal incompressible fluid flow}\label{s:ideal}

\subsection{The Cauchy invariants}\label{ss:cauchyinv}

The starting point of \citet{Cau27}
were the equations of incompressible fluid flow,
namely, the momentum equation and the incompressibility condition,
expressed in what is now called Eulerian coordinates.
In modern notation, they are, respectively:
\UB\partial_t\bv+(\bv\cdot\nabla)\bv=-\nabla p,\UE{EulerMom}
\UB\nabla\cdot\bv=0.\UE{Eulerincomp}
Here, $\bv$ is the flow velocity, $p$ pressure, $t$ time.

Cauchy transformed the equations by switching to Lagrangian coordinates,
denoted here by $\ba$. Let $\ba\mapsto\bx=\bx(\ba,t)$ be the Lagrangian map
from the initial to the current (i.e., at time $t$) Eulerian position, $\bx$,
of a fluid particle. The map satisfies
\UB\dot{\bx}=\bv(\bx,t),\qquad\bx(\ba,0)=\ba,\UE{xdotv}
where the dot denotes the Lagrangian time derivative. Since the left-hand
side of \rf{EulerMom} is the Lagrangian acceleration, \rf{EulerMom} implies
\UB\ddot{\bx}=-\nabla p=-\nL p((\nL\bx)^*)^{-1},\UE{LagDynamics}
where $\nL$ denotes the gradient in the Lagrangian variables and the star
denotes matrix transposition. Here, we have reexpressed the Eulerian pressure
gradient in terms of the Lagrangian pressure gradient, using the Jacobian
matrix $\nL\bx$ with entries $\nL_i x_j$. From here on, unless otherwise stated,
all space and time derivatives will be Lagrangian ones.

At this point, Cauchy took the Lagrangian curl of \rf{LagDynamics} and then
employed a {\it deus ex machina} to time-integrate the resulting equations.
Our approach involves only a minor modification, which eliminates some
of the magics: following \citet{We68}, we rewrite \rf{LagDynamics} as
\UB{\d\over\d t}\sum_{k=1}^3\dot{x}_k\nL x_k=\nL\left(
{1\over2}\,|\dot{\bx}|^2-p\right)\UE{equivDiffWeber}
and take the Lagrangian curl of \rf{equivDiffWeber} to find
\UB{\d\over\d t}\nL\times\sum_{k=1}^3\dot{x}_k\nL x_k
={\d\over\d t}\sum_{k=1}^3\nL\dot{x}_k\times\nL x_k=0.\UE{3members}
The expression in the middle is obtained by applying the Lagrangian
curl to each term in the sum and observing that no terms with second
spatial Lagrangian derivatives emerge. The importance of this will become clear
later. From \rf{3members} follows that the middle sum is a
time-independent vector. At $t=0$, by \rf{xdotv}, this vector is
seen to be the initial vorticity $\oz=\nL\times\vz$. Thus we obtain
Cauchy's equations in Lagrangian coordinates:
\be\sum_{k=1}^3\nL\dot{x}_k\times\nL x_k&=&\oz,\label{cauchy}\\
\det\,(\nL\bx)&=&1.\label{sdiff}\en
Here, \rf{sdiff} expresses the conservation of fluid volume.
The three components of the left-hand side of \rf{cauchy}, which are obviously
constant along any fluid particle trajectory, have been called
the {\it Cauchy invariants} \citep[see, e.g.,][]{YaZe}. Their conservation is
actually the three-dimensional generalisation of the two-dimensional
vorticity conservation.

As is well known, the unforced Euler equations \rf{Eulerincomp} are invariant
under Galilean transformations. More precisely, limiting ourselves
to solutions for which the velocity and the pressure are space-periodic,
if $(\bv(\bx,t),\,p(\bx,t))$ is a solution, then $(\bv(\bx-\bU t,t)+\bU,\break
p(\bx-\bU t,t))$ is also a solution to the same equation, for any uniform $\bU$.
This corresponds to changing the frame of reference to the one moving with the
velocity $\bU$. If, however, we try to describe the fluid
from an accelerating frame, we shall have to add a space-independent
force term. The latter is equivalent to adding to the pressure a term
linear in the spatial coordinates, but this is not consistent with
space periodicity. In contrast, when using Cauchy's Lagrangian equations
\rf{cauchy}--\rf{sdiff}, we can switch to an arbitrary accelerated frame,
provided it coincides with the original frame at $t=0$: this
just adds to the Lagrangian map a time-dependent uniform vector, which drops
out of \rf{cauchy}--\rf{sdiff}, because all the terms involve
gradients of the map. Of course, similar remarks can be made on the
Eulerian vorticity formulation of the Euler equation.

\subsection{The Taylor expansion in time}\label{ss:taylor}

The basic Lagrangian equations \rf{cauchy}--\rf{sdiff} do not give us
an explicit expression for the time derivative of either the Lagrangian map or of its gradient.
Nevertheless, they yield almost immediately explicit recurrence relations
for the time-Taylor coefficients of the Lagrangian map. This will be crucial
for the problems tackled here. For establishing these relations
it is convenient to work with the fluid particle
displacement field\break$\bxi(\ba,t)\equiv\bx(\ba,t)-\ba$.

In terms of the displacement, the Lagrangian equations \rf{cauchy} and
\rf{sdiff} take the form
\be&&\nL\times\dot{\bxi}+\sum_{k=1}^3\nL\dot{\xi}_k\times\nL\xi_k=\oz,
\label{3DDisplacement1}\\
&&\nL\cdot\bxi+\sum_{1\le i<j\le3}\left(\rule{0mm}{1em}
(\nL_i\xi_i)\nL_j\xi_j-(\nL_i\xi_j)\nL_j\xi_i\right)+\det(\nL\bxi)=0.
\label{3DDisplacement2}\en

Now, we assume that the displacement can be Taylor-expanded in the time
variable around $t=0$ for sufficiently short times:
\UB\bxi(\ba,t)=\sum_{s=1}^\infty\bxi^{(s)}(\ba)t^s.\UE{xiSeries}
Here the functions $\bxi^{(s)}$ are called the (time-)Taylor
coefficients. For the time being such an expansion is only formal.
The existence of all the time derivatives needed for the expansion will follow from the bounds
that are established in Sections~\ref{ss:L1} and \ref{ss:holder}.

Recurrence relations for the Taylor coefficients of the displacement
are obtained by substituting \rf{xiSeries} into \rf{3DDisplacement1} and
\rf{3DDisplacement2} and gathering all the terms containing a given
power $t^s$:
\be\nL\times\bxi^{(s)}&=&\oz\delta^s_1-\sum_{\subs{1\le k\le3,}{0<m<s}}
{m\over s}\,\nL\xi^{(m)}_k\times\nL\xi^{(s-m)}_k\nonumber\\
&=&\oz\delta^s_1-{1\over2}\sum_{\subs{1\le k\le3,}{0<m<s}}
{2m-s\over s}\,\nL\xi^{(m)}_k\times\nL\xi^{(s-m)}_k.\label{CauchyRecurrOpt}\\
\nL\cdot\bxi^{(s)}&=&\sum_{\subs{1\le i<j\le3}{0<m<s}}
\left(\rule{0mm}{1em}(\nL_j\xi^{(m)}_i)\nL_i\xi^{(s-m)}_j-(\nL_i\xi^{(m)}_i)
\nL_j\xi^{(s-m)}_j\right)\nonumber\\
&&-\sum_{\subs{i,j,k}{l+m+n=s}}
\varepsilon_{ijk}(\nL_i\xi^{(l)}_1)(\nL_j\xi^{(m)}_2)\nL_k\xi^{(n)}_3.
\label{sdiffRecurr}\en
Here $\varepsilon_{ijk}$ is the unit antisymmetric tensor and $\delta^s_{s'}$
is the Kronecker symbol; the second line of \rf{CauchyRecurrOpt} is obtained by grouping
the terms for the indices $m$ and $s-m$ in the sum in the right-hand side
of its first line. It is immediately seen from the two equations
\rf{CauchyRecurrOpt} and \rf{sdiffRecurr} for $s=1$ that $\bxi^{(1)}=\vz$.

A remarkable property of \rf{CauchyRecurrOpt}--\rf{sdiffRecurr}
is that their right-hand sides involve exclusively terms that are
quadratic and cubic in the gradients of the displacement and no higher
derivatives. This will be crucial in proving analyticity.

If we know all the Taylor coefficients up to order $s-1$, \rf{CauchyRecurrOpt}
and \rf{sdiffRecurr} give us the curl and the divergence of the next Taylor
coefficient of order $s$. A vector field, whose curl and divergence are known,
can be uniquely (for given boundary conditions) determined using
the Helmholtz--Hodge decomposition \citep[cf., e.g.,][]{AWe}. Applying it
to \rf{CauchyRecurrOpt}--\rf{sdiffRecurr}, we find
\UB\bxi^{(s)}=\nabla^{-2}(-\nL\times{\bf R}^{(1)}+\nL R^{(2)}),\UE{xis}
where $\nabla^{-2}$ denotes the inverse Laplacian in the Lagrangian coordinates
(taking into account the boundary conditions), and
${\bf R}^{(1)}$ and $R^{(2)}$ denote the right-hand sides of \rf{CauchyRecurrOpt}
and \rf{sdiffRecurr}, respectively.

\subsection{Proving analyticity: a case study using Fourier analysis}
\label{ss:L1}

To begin with, in Section~\ref{sss:elementary}, we shall prove the
time-analyticity of fluid particle trajectory by an elementary Fourier
method. We assume that the solution is $2\pi$-periodic in space and that
the spatial Fourier series of the initial vorticity
\UB\oz(\ba)=\sum_{\bf p}\widehat{\bom}_{\bf p}\,\e^{\i{\bf p}\cdot\ba},\UE{initom}
is absolutely convergent:
\UB\Gamma\equiv\sum_{\bf p}|\widehat{\bom}_{\bf p}|<\infty.\UE{absconv}
Here summation is over three-dimensional vectors $\bf p$ with integer
components. Because of the absolute convergence, clearly the
initial vorticity is continuous in $\ba$.

First, in Section~\ref{sss:elementary} we present a proof using elementary
transformations of the recurrence relations \rf{CauchyRecurrOpt}--\rf{sdiffRecurr}
in the spatial Fourier space. This ``elementary'' proof may,
however, leave the reader wondering why it is working. To prevent the trees
from hiding the forest, we then present, in Section~\ref{sss:abstract} another
slightly more transparent albeit more advanced ``normed-space'' proof in which the structure of the recurrence
relations is exploited to obtain bounds for the norms of the Taylor coefficients
in a suitable space of functions. The second
approach can and will be generalised to other normed function spaces
(initial vorticities that are analytic or H\"older-continuous), but it
does not yield estimates as sharp as the elementary proof.

\subsubsection{An elementary proof by transformations in Fourier space}
\label{sss:elementary}

Upon expanding the time-Taylor
coefficients of the displacement in spatial Fourier series
\UB\bxi^{(s)}(\ba)=\sum_{\bf p}\widehat{\bxi}^{(s)}_{\bf p}\,
\e^{\i{\bf p}\cdot\ba},\quad s=1,2,\ldots,\UE{TaylorFourier}
we can express the recurrence relations in the Fourier representation.
Eq.~\rf{CauchyRecurrOpt} becomes then
\UB{\bf p}\times\widehat{\bxi}^{(s)}_{\bf p}=
-\i\widehat{\bom}_{\bf p}\delta^s_1-{\i\over2}\sum_{0<m<s,\ \bf r}{2m-s\over s}\,
(\widehat{\bxi}^{(m)}_{\bf r}\cdot\widehat{\bxi}^{(s-m)}_{{\bf p}-{\bf r}})
\,({\bf r}\times({\bf p}-{\bf r})).\UE{CauchyTF}
Similarly, it is elementary but a bit tedious (see Appendix \ref{appA})
to rewrite \rf{sdiffRecurr} as
\be{\bf p}\cdot\widehat{\bxi}^{(s)}_{\bf p}
&=&{\i\over2}\sum_{0<m<s,\ \bf r}({\bf r}\times({\bf p}-{\bf r}))
\cdot\left(\widehat{\bxi}^{(m)}_{\bf r}\times\widehat{\bxi}^{(s-m)}_{{\bf p}-{\bf r}}\right)\nonumber\\
&&+{1\over6}\sum_{\subs{n_1+n_2+n_3=s}{{\bf r}_1+{\bf r}_2+{\bf r}_3=\bf p}}
[{\bf r}_1,{\bf r}_2,{\bf r}_3]\,\left[\widehat{\bxi}^{(n_1)}_{{\bf r}_1},
\widehat{\bxi}^{(n_2)}_{{\bf r}_2},\widehat{\bxi}^{(n_3)}_{{\bf r}_3}\right],\label{sdiffTF}\en
where $[\cdot,\cdot,\cdot]$ denotes the triple product.

For ${\bf p}=0$ equations \rf{CauchyTF} and \rf{sdiffTF} are trivially
satisfied. The mean (over the periodicity cell) of the velocity is
easily seen to be conserved. Without loss of generality, we henceforth assume
that this mean vanishes initially. Consequently, the mean velocity and the
mean displacement will also vanish at any time, and thus
$\widehat{\bxi}^{(s)}_0=0$ for all~$s$.

For ${\bf p}\ne0$ the coefficient $\widehat{\bxi}^{(s)}_{\bf p}$ can be
uniquely determined, as we know its components normal and parallel to
the wave vector $\bf p$, by virtue of \rf{CauchyTF} and \rf{sdiffTF},
respectively. (In the Fourier representation, the Helmholtz--Hodge
decomposition is indeed just a geometric decomposition of a vector
into transverse and longitudinal components, perpendicular and parallel
to the wave vector $\bf p$.)
Let $\widetilde{\bf R}_{\bf p}^\perp$ and $\widetilde R_{\bf p}^\parallel$
denote the right-hand side of \rf{CauchyTF} and \rf{sdiffTF}, respectively. Then
\UB\widehat{\bxi}^{(s)}_{\bf p}=|{\bf p}|^{-2}\left(-{\bf p}\times
\widetilde{\bf R}_{\bf p}^\perp+\widetilde R_{\bf p}^\parallel\,{\bf p}\right).\UE{xip}

We are now using recurrence relations \rf{CauchyTF} and \rf{sdiffTF} to derive
bounds for the Fourier coefficients of the time-Taylor series. Eq.~\rf{CauchyTF}
and the observation that $|2m -s|/s \le 1$ for all $0<m<s$ imply a bound
for the transverse component:
\UB|{\bf p}\times\widehat{\bxi}^{(s)}_{\bf p}|\le
|\widehat{\bom}_{\bf p}|\delta^s_1+{1\over2}\sum_{0<m<s,\ \bf r}|{\bf r}|\,|{\bf p}-{\bf r}|\,
|\widehat{\bxi}^{(m)}_{\bf r}\cdot\widehat{\bxi}^{(s-m)}_{{\bf p}-{\bf r}}|.\UE{FCauchyTFbound}
A bound for the longitudinal component follows from \rf{sdiffTF}:
\be|{\bf p}\cdot\widehat{\bxi}^{(s)}_{\bf p}|&\le&{1\over2}\sum_{m,\bf r}
|{\bf r}|\,|{\bf p}-{\bf r}|\,|\widehat{\bxi}^{(m)}_{\bf r}\times
\widehat{\bxi}^{(s-m)}_{{\bf p}-{\bf r}}|\nonumber\\
&&+{1\over6}\sum_{\subs{n_1+n_2+n_3=s}{{\bf r}_1+{\bf r}_2+{\bf r}_3=\bf p}}
|{\bf r}_1|\,|{\bf r}_2|\,|{\bf r}_3|\,|\widehat{\bxi}^{(n_1)}_{{\bf r}_1}|\,
|\widehat{\bxi}^{(n_2)}_{{\bf r}_2}|\,
|\widehat{\bxi}^{(n_3)}_{{\bf r}_3}|.\label{sdiffTFbound}\en

We now introduce the generating function
\UB\zeta(t)\equiv\sum_{s=1}^\infty\left(\sum_{\bf p}|{\bf p}|\,|\widehat{\bxi}^{(s)}_{\bf p}|\right)t^s.\UE{genfu}
 From \rf{genfu} and the bounds \rf{FCauchyTFbound} and \rf{sdiffTFbound}
for the transverse and longitudinal components,
\be\zeta(t)&\le&\sum_{s=1}^\infty\sum_{\bf p}\left(|{\bf p}\times\widehat{\bxi}^{(s)}_{\bf p}|
+|{\bf p}\cdot\widehat{\bxi}^{(s)}_{\bf p}|\right)t^s\nonumber\\
&\le&\sum_{\bf p}|\widehat{\bom}_{\bf p}|t+{1\over\sqrt2}\sum_{s=1}^\infty
\sum_{\subs{0<m<s,}{\bf r}}|{\bf r}|\,|{\bf p}-{\bf r}|
|\widehat{\bxi}^{(m)}_{\bf r}|\,|\widehat{\bxi}^{(s-m)}_{{\bf p}-{\bf r}}|\,
t^s\nonumber\\
&&+{1\over6}\sum_{s=1}^\infty\,\sum_{\subs{n_1+n_2+n_3=s}{{\bf r}_1+{\bf r}_2+{\bf r}_3=\bf p}}
|{\bf r}_1|\,|{\bf r}_2|\,|{\bf r}_3|\,|\widehat{\bxi}^{(n_1)}_{{\bf r}_1}|\,
|\widehat{\bxi}^{(n_2)}_{{\bf r}_2}|\,
|\widehat{\bxi}^{(n_3)}_{{\bf r}_3}|\,t^s\nonumber\\
&=&T+{\zeta^2(t)\over\sqrt 2}+{\zeta^3(t)\over6},\label{final}\en
where $T\equiv\Gamma t$ is a dimensionless time; we have used
the geometric interpretation of scalar and vector products and applied
the inequalities $1\le|\cos\theta|+|\sin\theta|\le\sqrt 2$ that hold
for any angle $\theta$.

\begin{figure}
\centerline{\includegraphics[width=\textwidth,height=3in]{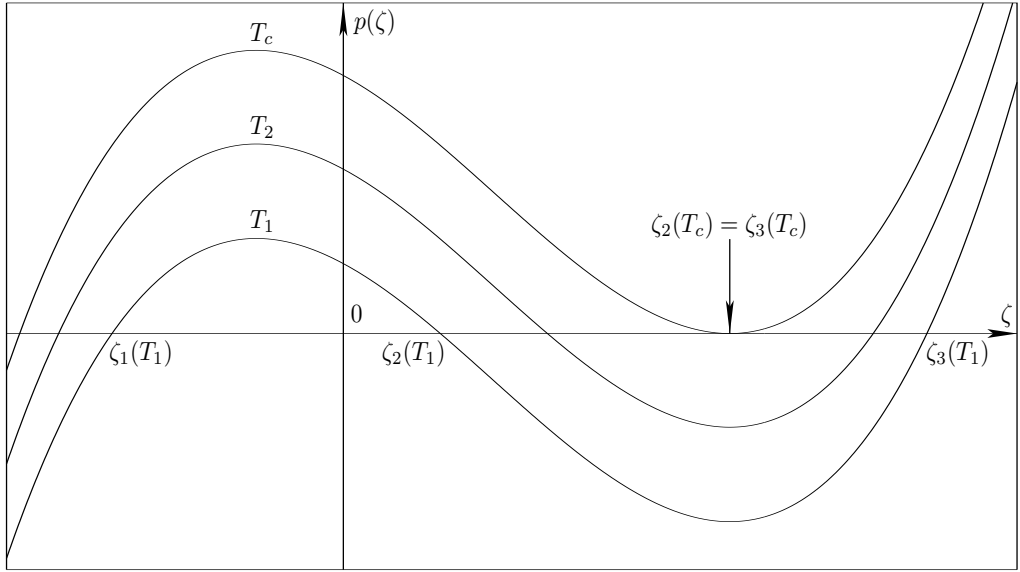}}
\caption{A sketch of graphs of the polynomial $p(\zeta)$ \rf{poly} for three
values of $T$, such that $T_1<T_2<T_c$, illustrating the behaviour of real roots of $p$
(the roots are shown as the points of intersection of the graph of $p(\zeta)$
and the horizontal axis). On increasing $T$, the graph slides up as a rigid
curve. As a result, the roots $\zeta_1$ and $\zeta_3$ move to the left (i.e.,
become smaller), and $\zeta_2$ moves to the right (i.e., becomes larger).
At the critical value $T_c$, for which $\partial p/\partial\zeta=0$ and
the discriminant $\Delta$ vanishes, the two roots $\zeta_2$ and $\zeta_3$
collide and disappear (with emergence of a pair of complex conjugate roots,
which cannot be shown in the plot).}
\label{fig1}\end{figure}

The inequality \rf{final} is equivalent to
\UB p(\zeta)\equiv{\zeta^3\over6}+{\zeta^2\over\sqrt 2}-\zeta+T\ge0.\UE{poly}
The discriminant of the polynomial $p$ \rf{poly},
\UB\Delta=-{3\over4}\,T^2-{5\over\sqrt{2}}\,T+{7\over6},\UE{discriminant}
is positive at small times, and thus $p$ has three real roots $\zeta_i$.
By Vi\`ete's theorem, the product of the roots is negative;
since the extrema of the graph of $p$ lie at points of different signs,
we infer that $\zeta_1<0<\zeta_2<\zeta_3$. Differentiating in $T$
the equation for the roots of $p$, we find
\UB{\partial\zeta_i\over\partial T}=-\left(\left.{\partial p\over
\partial\zeta}\right|_{\zeta=\zeta_i}\right)^{-1}.\UE{dpdT}
Consequently, on increasing $T$, the root $\zeta_2$ monotonically
increases and $\zeta_3$ monotonically decreases till the two roots collide
(i.e., $\zeta_2=\zeta_3$) giving birth to a pair of complex roots
(see Fig.~\ref{fig1}). This occurs at a critical value
\UB T=T_c=(8-5\sqrt{2})/3\approx0.3096\UE{Tcrit}
for which $\Delta=0$.
For $T<T_c$, inequality \rf{poly} implies $\zeta\le\zeta_2$. Therefore,
the series $\sum_{s,\bf p}|\widehat{\bxi}^{(s)}_{\bf p}(t)|t^s$
converges for $t<t_c\equiv T_c/\Gamma$, and hence at any Lagrangian
position $\bf a$ the displacement
$\bxi(\ba,t)=\sum_{s,\bf p}\widehat{\bxi}^{(s)}_{\bf p}(t)\,\e^{\i{\bf p}\cdot\ba}\,t^s$
is analytic in complex time $t$ in the disk $|t|<t_c$.

We have just derived the bound $\zeta(t)\le\zeta_2$ under the assumption
that the displacement is given by its time-Taylor series. The assumption has
enabled us to work with the infinite series, but this is actually not required.
A complete mathematical proof proceeds as follows. First, by essentially
repeating the arguments of this section, one can show that for $t<t_c$
any partial sum of the series \rf{genfu} is also bounded by $\zeta_2$. This
bound being uniform, we conclude that the series \rf{xiSeries} for $\bxi(\ba,t)$
converges in the disk $|t|<t_c$ for any $\ba$, and defines an analytic function
of time. Consequently, the $s$-th coefficient of the generating function
\rf{genfu} is bounded by $\zeta_2t_c^{-s}$. Second, one substitutes
the truncated series \rf{xiSeries} into equations \rf{cauchy} and \rf{sdiff}.
Using our bound for the coefficients in \rf{xiSeries}, we easily find that
for $t<t_c$ the discrepancies decay as $CS^2(t/t_c)^S(1-t/t_c)^{-1}$ and
$CS^3(t/t_c)^S(1-t/t_c)^{-1}$, respectively,
where $S$ denotes the number of the terms kept and $C$ is a constant. Thus,
in the disk of convergence, the {\it formal solution} that we have constructed
in this section is a classical solution to Cauchy's equations.

\subsubsection{A normed-space proof}\label{sss:abstract}

It may be instructive to give an alternative proof, that is not tailored
to the specific functional space to which the initial vorticity belongs,
while emphasizing the key ingredients of such a proof.

It is crucial for further developments to reexpress the pair of recurrence
relations \rf{CauchyRecurrOpt}--\rf{sdiffRecurr} as a single recurrence
for the (Lagrangian) gradients of the Taylor coefficients (the latter are
vectors, so their gradients are 2-tensors). We take the curl
of \rf{CauchyRecurrOpt} and the gradient of \rf{sdiffRecurr}, solve
the resultant Poisson equation and finally take one (tensor) gradient to obtain
\be\nL_\mu\xi_\nu^{(s)}&=&\nL_\mu v_\nu^{\rm init}\delta^s_1+\sum_{\subs{1\le j\le3,}{j\ne\nu}}
{\cal C}_{\mu j}\left(\sum_{\subs{1\le k\le3,}{0<m<s}}{2m-s\over s}\,
(\nL_\nu\xi^{(m)}_k)\,\nL_j\xi^{(s-m)}_k\right)\nonumber\\
&&+\,{\cal C}_{\mu\nu}\left(\sum_{\subs{1\le i<j\le3}{0<m<s}}
\left(\rule{0mm}{1em}(\nL_j\xi^{(m)}_i)\nL_i\xi^{(s-m)}_j-(\nL_i\xi^{(m)}_i)
\nL_j\xi^{(s-m)}_j\right)\right.\nonumber\\
&&\left.-\sum_{\subs{i,j,k}{l+m+n=s}}
\varepsilon_{ijk}(\nL_i\xi^{(l)}_1)(\nL_j\xi^{(m)}_2)\nL_k\xi^{(n)}_3\right).
\label{bigXis}\en
Here, ${\cal C}_{ij}\equiv\nabla^{-2}\nL_i\nL_j$ and $\nL_i\nL_j$ denotes
the second-order partial derivative $\partial^2/\partial a_i\partial a_j$.
The operator $\nabla^{-2}$, as already stated, is the inverse of the (Lagrangian)
Laplacian: given a periodic function $f$ with zero spatial mean over the periodic
cell, $\nabla^{-2}f$ is defined as the unique periodic function $\psi$ with
zero spatial mean, solving $\nabla^2\psi=f$.
The structure of \rf{bigXis} has an important property: for a given order $s$,
the 9-component tensor $\nL\bxi^{(s)}$ is expressed in terms of products
of gradients of lower-order Taylor coefficients, and of the operators
${\cal C}_{ij}$ acting on these products.

Operators such as ${\cal C}_{ij}\equiv\nabla^{-2}\nL_i\nL_j$ have been
used in potential theory and applied mathematics since at least the
beginning of the 20th century (see Section~\ref{ss:holder}). They are
instances of what is now sometimes called Calderon--Zygmund
operators. In view of the importance of these operators here, we shall
call them the ``fundamental Calderon--Zygmund'' (fCZ) operators.

Suppose a functional norm has the following two
properties: $(i)$ the norms of the fCZ operators are bounded; $(ii)$
the norm of a product is bounded by the product of norms of its factors. Then
clearly relation \rf{bigXis} enables us to control the norm
of the gradient of each subsequent Taylor coefficient. Actually, we have already
used, albeit without stressing it, this structural property of the recurrence
relations for the gradients of Taylor coefficients in the elementary proof.

We now define a norm which, as we shall show, has the two required properties:
For a scalar-valued space-periodic function
\UB f(\ba)=\sum_{\bf p}\widehat f_{\bf p}\,\e^{\i{\bf p}\cdot\ba}\UE{fFourier}
we set
\UB\|f\|\equiv\sum_{\bf p}|\widehat f_{\bf p}|.\UE{Snorm}
Following a standard approach, we then set
\UB\|{\bf f}\|\equiv\max_j\|f_j\|\UE{Vnorm}
for a space-periodic vector field $\bf f$, and
\UB\|\nL{\bf f}\|\equiv\max_i\|\nL_i{\bf f}\|=\max_{i,j}\|\nL_if_j\|\UE{Tnorm}
for the tensor $\nL\bf f$. Actually, \rf{Tnorm} is equivalent to the norm
$\sum_{\bf p}|{\bf p}|\,|\widehat{\bxi}^{(s)}_{\bf p}|$ of $\nL\bxi^{(s)}$,
used in the elementary proof in the previous section. Other norms having
the two required properties will be used in subsequent sections.

$i$. {\it Boundedness of fCZ operators ${\cal C}_{ij}$}

The operator ${\cal C}_{ij}$ is a composition of two
derivatives and the inverse Laplacian in Lagrangian variables. For
given boundary conditions, applying the inverse Laplacian amounts to
solving a Poisson equation. This is a nonlocal operation. Formally
counting the number of differentiations involved in the fCZ
operator, we see that it is a kind of zeroth-order nonlocal
differential operator --- in the theory of partial differential
operators such an operator is called a pseudodifferential operator
of order zero. For the norm $\|\cdot\|$, the fCZ operator
${\cal C}_{ij}$ is bounded; actually its norm is unity. Indeed, for any
space-periodic function $f$ expanded in a Fourier series \rf{fFourier}, we have
\UB\|{\cal C}_{ij}f\|=\left\|\sum_{{\bf p}\ne0}{p_ip_j\over|{\bf p}|^2}
\widehat f_{\bf p}\,\e^{\i{\bf p}\cdot\ba}\right\|
=\sum_{{\bf p}\ne0}{|p_ip_j|\over|{\bf p}|^2}|\widehat f_{\bf p}|
\le\sum_{\bf p}|\widehat f_{\bf p}|=\|f\|,\UE{boundedOperator}
as stated.

$ii$. {\it A bound for the norms of a product}

We now show that the norm \rf{Snorm} has the second desired
property, i.e. the norm of a product is bounded by the product of norms
of the factors. Consider two functions $f$ and $g$, whose Fourier coefficients
are $\widehat f_{\bf p}$ and $\widehat g_{\bf p}$, respectively. Clearly,
\UB\|fg\|=\sum_{\bf p}\left|\sum_{\bf n}\widehat f_{\bf n}\widehat g_{\bf p-n}\right|\le
\sum_{\bf p}\sum_{\bf n}|\widehat f_{\bf n}|\,|\widehat g_{\bf p}|=\|f\|\|g\|,\UE{prodBound}
as required. A linear functional space equipped with a norm with this property
is called an algebra.

The right-hand side of \rf{bigXis} involves three sums. For a given index
of summation $m$, each of the first two sums involve 6 quadratic products
of various $\nL_{i'}\xi^{(j')}$; for given $l$ and $m$, the third sum involves
6 cubic products of such components. Using this and applying $i$ and $ii$, we bound
the norm of the various terms in the right-hand side of \rf{bigXis} and obtain
\be\|\nL\bxi^{(s)}\|&\le&\|\nL\vz\|\delta_1^s+12\sum_{m+n=s}\|\nL\xi^{(m)}\|\,\|\nL\xi^{(n)}\|\nonumber\\
&&+\,6\sum_{l+m+n=s}\|\nL\xi^{(l)}\|\,\|\nL\xi^{(m)}\|\,\|\nL\xi^{(n)}\|.\label{finalXiBound}\en
We then introduce the generating function
\UB\widetilde\zeta(t)\equiv\sum_{s=1}^\infty\|\nL{\bxi}^{(s)}\|t^s,\UE{genfunc}
similarly to what we did in the elementary proof. Multiplying
\rf{finalXiBound} by $t^s$ and summing over $s\ge1$ yields
\UB6\widetilde\zeta^3+12\widetilde\zeta^2-\widetilde\zeta+\|\nL\vz\|t\ge0\,.\UE{abstrPoly}
Inequality \rf{abstrPoly} has the same nature as \rf{poly}, and the remainder
of the proof is along the same lines as in the elementary proof.
Time analyticity is guaranteed by \rf{abstrPoly} until the critical time
$\approx 0.0204/\|\nL\vz\|$.

Neglecting in the normed-space proof fine details of equations arising
in the Fourier space in order to gain a higher generality has resulted
in a quantitative (but not qualitative) deterioration of the bounds (for instance,
the coefficient in the cubic term in \rf{abstrPoly} is 36 times larger
than the one in \rf{poly}). Although the dimensionless time is defined
differently in the ``elementary'' and the ``normed-space'' proofs, one easily
finds, using the inequality $\|\nL{\bxi}^{(s)}\|\le\Gamma\le6\|\nL{\bxi}^{(s)}\|$,
that the bound for the radius of convergence of the Taylor series is at
least 2.53 times larger in the ``elementary'' proof.

\subsection{Analytic initial vorticity}\label{ss:analytic}

Considering again $2\pi$-periodic flows, we now assume
that the initial condition is an analytic function of the three
space variables. We shall prove that the Lagrangian map is a space-time
analytic function over the whole periodic domain and until times that again
depend only on the initial velocity gradient.

For this, we employ a Gevrey-type norm,
that takes advantage of the exponential decrease with wavenumber
of the spatial Fourier coefficients of analytic functions. For $\sigma>0$ and
a scalar-valued function $f$ \rf{fFourier}, we define
\UB\|f\|_\sigma\equiv\sum_{\bf p}\e^{\sigma|{\bf p}|}
\,|\widehat f_{\bf p}|,\UE{GevreyNorm}
and for a vector field $\bf f$,
\UB\|\nL{\bf f}\|_\sigma\equiv\max_{i,j}\|\nL_if_j\|_\sigma.\UE{analNorm}

Since the initial velocity $\vz(\ba)$ is a periodic analytic function,
for some $\sigma>0$,
\UB\|\nL\vz\|_\sigma=\max_{ij}\sum_{\bf p}
|p_i|\e^{\sigma|{\bf p}|}\,|(\widehat\bv_{\bf p})_j|\UE{nOm}
is finite. This guarantees that the Fourier series for the initial velocity and
\rf{initom} for the initial vorticity converge for complex $\ba$ in the strip
$|{\rm Im}\,z|<\sigma$ around the real domain. Moreover, the supremum
of $\sigma$, for which the sum \rf{nOm} converges, can be identified as
the distance from the real space $\R^3$ to the nearest singularity
of the analytic function $\nL\vz(\ba)$ in the complex three-dimensional space.
As in Section~\ref{sss:abstract}, we consider the generating function
\UB\widetilde\zeta_\sigma(t)
\equiv\sum_{s=1}^\infty\|\nL{\bxi}^{(s)}\|_\sigma\,t^s.\UE{sigmagenfunc}
If we show that the generating function is finite until $t=t_c(\sigma)$,
spatio-temporal analyticity
until that time follows. This will hold if we demonstrate the two required
properties of the norm $\|\cdot\|_\sigma$. Their proofs are similar
\textit{mutatis mutandis} to those
given in the previous subsection in \rf{boundedOperator} and \rf{prodBound}.
% demonstrations are similar
%to the ones presented in Section~\ref{sss:abstract}.
%
%$i$. {\it Boundedness of the fCZ operators ${\cal C}_{ij}$}
%
%For any analytic space-periodic function $f$ with Fourier coefficients
%$\widehat f_{\bf p}$, we have
%\UB\|{\cal C}_{ij}f\|_\sigma=\left\|\sum_{{\bf p}\ne0}{p_ip_j\over|{\bf p}|^2}
%\widehat f_{\bf p}\,\e^{\i{\bf p}\cdot\ba}\right\|_\sigma
%=\sum_{{\bf p}\ne0}{|p_ip_j|\over|{\bf p}|^2}\,\e^{\sigma|{\bf p}|}\,|\widehat f_{\bf p}|
%\le\sum_{\bf p}\e^{\sigma|{\bf p}|}|\,\widehat f_{\bf p}|=\|f\|_\sigma.\UE{sigmaBoundedOperator}
%
%$ii$. {\it Bounds for products}
%
%By virtue of the triangle inequality, for any two analytic space-periodic
%functions $f$ and $g$ with Fourier coefficients $\widehat f_{\bf p}$ and
%$\widehat g_{\bf p}$, respectively,
%\UB\|fg\|_\sigma=\sum_{\bf p}\left|\sum_{\bf n}\widehat f_{\bf n}\widehat g_{\bf p-n}\right|
%\e^{\sigma|{\bf p}|}\le\sum_{\bf p}\sum_{\bf n}|\widehat f_{\bf n}|\,|\widehat g_{\bf p-n}|
%\e^{\sigma(|{\bf n}|+|{\bf p-n}|)}=\|f\|_\sigma\|g\|_\sigma.\UE{sigmaProdBound}

Then, by repeating the normed-space proof almost literally, we derive the inequality
\UB6\widetilde\zeta_\sigma^3+12\widetilde\zeta_\sigma^2-\widetilde\zeta_\sigma+t\,\|\nL\vz\|_\sigma\ge0,\UE{analytPoly}
from which analyticity in space and time follows.

More careful estimations, as in the elementary proof
of Section~\ref{sss:elementary} show that the generating function
\UB\zeta_\sigma(t)\equiv\sum_{s=1}^\infty\left(\sum_{\bf p}|{\bf p}|\,\e^{\sigma|{\bf p}|}\,
|\widehat{\bxi}^{(s)}_{\bf p}|\right)t^s\UE{analGenfu}
satisfies the inequality
\UB{\zeta_\sigma^3\over6}+{\zeta_\sigma^2\over\sqrt 2}-\zeta_\sigma+t\,\Gamma_\sigma\ge0,\UE{sigmapoly}
where $\Gamma_\sigma\equiv\|\oz\|_\sigma$.
Reasoning precisely as in Section~\ref{sss:elementary}, we infer that
$\zeta_\sigma$ remains bounded and therefore the Lagrangian map is analytic
for $t<T_c/\Gamma_\sigma$ where $T_c \approx 0.3096$, as given by \rf{Tcrit},
and its singularities, if any, are at a distance at least $\sigma$
from the real space. Letting $\sigma\searrow0$, and observing that
$\Gamma_\sigma$ is an increasing
continuous function of $\sigma$, we conclude that $\bxi(\ba,t)$ is guaranteed
to possess spatio-temporal analyticity for $t<t_c\equiv T_c/\Gamma_0$.

For analytic initial data, analyticity of the Lagrangian map in space and time
is not surprising, given that in Eulerian coordinates it has been
known for a long time that, with analytic initial data, the solution is
analytic in space and time (\citealp{Be76}; spatial analyticity of space-periodic
solutions was further considered by \citealp{LeOl}, \citealp{KV09} and
\citealp{Zh11}).
The Lagrangian map is the solution to $\dot{\bx}=\bv(\bx,t)$ with $\bx(0)=\ba$,
where $\bv(\bx,t)$ is the Eulerian velocity. From this one can in principle
establish analyticity of the Lagrangian map, using the theory
of ODEs of complex variables.

Nevertheless, with our derivation, we have gained something essential:
the latter proof, by way of Eulerian analyticity, will not work beyond complex
times whose moduli exceed the radius of convergence of the Taylor
series in $t$ for any Eulerian $\bx$. This Eulerian radius of convergence
does not depend only on the {\it initial velocity gradient}, as is the case
of the time of guaranteed analyticity by our Cauchy-formulation-based technique,
but also on the {\it initial velocity itself}. This can be seen by adding
to the initial flow a large uniform velocity $\bU$. As a consequence,
a singularity at a given complex space location $\bx_\star$ within a distance
$\delta$ of the real domain will be swept past a fixed (Eulerian) observer
at speed $\bU$ and, to leading order, it will also be seen as a complex
singularity in time at a location whose imaginary part is $\delta/|\bU|$.
Thus as $|\bU|$ is increased, the Eulerian radius of convergence decreases.
Note that the argument is basically equivalent to the use of the Taylor
hypothesis, discussed in the Introduction.

This argument has an important implication for numerical simulation of
ideal incompressible flow, particularly when high precision is needed,
e.g. to study the appearance of possible singularities. The simplest
way to achieve high precision in temporal schemes is to use a
high-order scheme, for example a high-order truncated Taylor series in
time, the time-step of which cannot normally exceed the radius of
convergence of the Taylor series. In Eulerian coordinates this may
result in a much smaller permissible time step than in Lagrangian
coordinates. We shall discuss this issue elsewhere.

Finally, let us remark that the use of Lagrangian coordinates has proved useful
for establishing the spatial analyticity of solutions in a 3D fluid region
with a boundary, and deriving bounds for the rate of decay of their radius
of spatial analyticity \citep{KV11}.

\subsection{H\"older-continuous initial vorticity}\label{ss:holder}

Again, we here assume spatial periodicity
and denote the elementary periodicity cell by $\T^3\equiv[0,\,2\pi]^3$.
We assume that the initial vorticity is H\"older-continuous with
some exponent $0<\alpha<1$. This means that its H\"older norm
\UB|\oz|_\alpha\equiv\max_{\ba\in\T^3}\,|\oz(\ba)|+\sup_{\ba_1,\,\ba_2\in\T^3,\,\ba_1\ne\ba_2}
{|\oz(\ba_1)-\oz(\ba_2)|\over|\ba_1-\ba_2|^\alpha}\UE{nullAlphaNorm}
is finite. Such a function is continuous but may be non-differentiable,
and its increments between two neighbouring points $\ba_1$ and $\ba_2$ may
behave as $|\ba_1-\ba_2 |^\alpha$.

For such flows we give a proof of the finite-time temporal
analyticity of the Lagrangian map. It is similar to the proof by
\citet{FriZh}, but builds on the material of Section~\ref{sss:abstract} and
is more detailed, in order to have a self-contained presentation for
fluid dynamicists.

$i$. {\it Boundedness of the fCZ operators ${\cal C}_{ij}$}

Here we consider the following statement: {\it Let a function $f(\ba)$
be periodic and H\"older-continuous of exponent $0<\alpha<1$, then}
\UB|\nabla^{-2}\nL_i\nL_j f(\ba)|_\alpha\le\Theta_\alpha|f|_\alpha,\UE{HolCijBou}
{\it where $\Theta_\alpha$ does not depend on $f$ and is bounded by}
$C/\alpha$. In the language of electrostatics, an equivalent statement would be
that if the charge density is H\"older-continuous, so is the gradient
of the electric field. In hydrodynamics, an analogous result is that, if
the vorticity is H\"older-continuous, so is the gradient of the velocity,
the latter being given by a Biot-Savart integral. Note that a vorticity
that is merely continuous can give rise to an arbitrarily large velocity
gradient. Hence such a
vorticity could be shear-stretched arbitrarily fast, leading to a
vanishing interval of assured regularity.

Closely related results have been proved by mathematicians in the early
twentieth century. \citet{Li25} stated a similar result in the entire space
$\R^3$, giving a reference to \citet{Ko} and mentioning that
``it does not present serious difficulties but requires, as is known,
a series of rather cumbersome thoughts.'' He used it to prove various results
regarding the initial-value problem for the 3D incompressible Euler
equations, including its well-posedeness for a finite time when the initial
vorticity is H\"older-continuous \citep{Li27}. \citet{Ch92} used the boundedness
of the fCZ operators (essentially, in the Eulerian framework) to prove
that when the initial vorticity is H\"older-continuous, the fluid particle
trajectories are infinitely differentiable in time.

There exist two types of proofs of the boundedness of the fCZ
operators similar to \rf{HolCijBou}. When the problem is set in the
entire $\R^3$ space, the Green function for the Poisson equation is
basically the Coulomb potential $1/r$ and proofs can be given by
bounding suitable singular integral operators in physical space. With
periodicity conditions, the problem is somewhat more cumbersome,
although the Green function still behaves as $1/r$ at distances $r$
small compared to the period. The method, using so-called Poisson integrals,
\citep[see, e.g.,][Chapters 3 and 4]{Ste} can handle both instances;
it is rather elementary, but definitely too long
for presenting it here, even in a summarised form. To give the reader
an idea of why the boundedness holds, we shall discuss
the case of the entire space and for simplicity assume that the given
function $f$ is not only H\"older-continuous of exponent $\alpha$, but
also vanishes outside a bounded domain (i.e., is compactly supported).
We here follow \citet[Chapter 4]{GT98} and \citet[Chapter 4]{MB02}.

The solution to the three-dimensional Poisson equation
\UB\nabla^2\phi=f\UE{Poisson}
is a convolution
\UB\phi(\ba)=\int f(\by)G(\ba-\by)\,\d\by,\UE{convo}
where the integral is taken over the whole space and the Green function $G(\ba)$ is
\UB G(\ba)=-{1\over4\pi|\ba|}.\UE{spaceGreen}

In spite of the singular behaviour of the Green function at $\ba=0$,
the integral \rf{convo} converges. Furthermore, we can differentiate
it once with respect to $\ba$ by just differentiating the Green function to obtain
\UB\nL_i\phi(\ba)=\int f(\by)\,\nL_iG(\ba-\by)\,\d\by.\UE{space1der}
If however we wish to differentiate twice, we cannot proceed in the
same way, since the second derivatives of the Green function have a strong,
$|\ba|^{-3}$, singularity which, in three dimensions, could lead to a
logarithmic divergence at short distance. This is overcome in the following
way. In principle to differentiate a convolution integral such as \rf{space1der},
we can differentiate either of the factors. So, we can apply one of the
differentiations to $f$ rather than to the Green function. For this
we need to assume temporarily that $f$ is differentiable. As we shall see, the
final expression for the second derivative of $\phi$ does not involve
derivatives of $f$ and holds for functions which are merely H\"older-continuous.
Standard arguments allow us then to justify the formula for such functions. Thus,
\UB\nL_i\nL_j\phi(\ba)=\int{\partial\over\partial a_j}f(\ba-\by')\,
{\partial\over\partial y'_i}\,G(\by')\,\d\by',\UE{space2der}
where we have
changed variables from $\by$ to $\by'\equiv\ba-\by$. We note that
\UB{\partial\over\partial a_j}f(\ba-\by')=
-{\partial\over\partial y'_j}f(\ba-\by'),\UE{nonlabel}
and subtract a {\it counterterm}, which is clearly zero, from \rf{space2der}
to get
\UB\nL_i\nL_j\phi(\ba)=-\int{\partial\over\partial y'_j}
\,\Big(f(\ba-\by')-f(\ba)\Big)\,
{\partial\over\partial y'_i}\,G(\by')\,\d\by'.\UE{spaceder}
Now, we return to the original spatial variable $\by$
and observe that, since $f$ was assumed to vanish outside of a bounded
set, called its support, we do not have to integrate over the whole
space. It is convenient to chose a set $\Omega$ with a smooth boundary,
for example a sphere, whose boundary $\partial \Omega$ is at a
strictly positive distance from the support of $f$. Using the
Gauss--Ostrogradsky divergence theorem, we obtain
\UB\nL_i\nL_j\phi=\!\int_\Omega\!\Big(\!f(\by)-f(\ba)\Big)\nL_i\nL_jG(\ba-\by)\,\d\by
-f(\ba)\!\int_{\partial\Omega}\nL_iG(\ba-\by)\nu_j(\by)\,\d\sigma(\by),
\UE{derivative}
where $\bm\nu$ is the outward unit vector normal to $\partial\Omega$
(note that $f(\by)=0$ for $\by$ on $\partial\Omega$).

The first integral in the right-hand side is nearly what we would have obtained
by illegitimate naive double differentiation of the Green function. However,
$f$ appears now only through its increments. By H\"older continuity, these are
bounded by $|f|_\alpha |\by-\ba|^\alpha$, thus preventing the aforementioned
logarithmic divergence. Specifically, we have
\UB\left|\int_\Omega\!(f(\by)-f(\ba))\nL_i\nL_jG(\ba-\by)\,\d\by\right|
\le c|f|_\alpha\int_\Omega|\by-\ba|^{\alpha-3}\,\d\by\le{C\over\alpha}\,|f|_\alpha,
\UE{intbou}
where $c$ and $C$ are suitable positive constants. To bound the
H\"older norm of the first integral in \rf{derivative}, we also need
to bound its increments. Establishing such bounds is a bit technical but quite
straightforward \citep[see, e.g.,][Chapter 4]{GT98}.

As to the second, surface integral in \rf{derivative} --- since the boundary
$\partial\Omega$ is at a positive distance from the set where $f\ne 0$, this
integral is an indefinitely differentiable function on the support of $f$.
We have thus established the inequality \rf{HolCijBou}, albeit
in the framework of an entire-space formulation
of the Poisson problem with a right-hand side that has bounded support.

$ii$. {\it Bounds for products of H\"older-continuous functions}

The other property we need is straightforward. Let $f$ and $g$ be
H\"older-continuous functions of exponents $\alpha$. Then
\be|fg|_\alpha&\le&\max_{\ba\in\T^3}\,|f(\ba)|\,\max_{\ba\in\T^3}\,|g(\ba)|
+\max_{\ba_1\in\T^3}\,|g(\ba_1)|\sup_{\ba_1,\,\ba_2\in\T^3,\,\ba_1\ne\ba_2}
{|f(\ba_1)-f(\ba_2)|\over|\ba_1-\ba_2|}\nonumber\\
&&+\max_{\ba_2\in\T^3}\,|f(\ba_2)|\sup_{\ba_1,\,\ba_2\in\T^3,\,\ba_1\ne\ba_2}
{|g(\ba_1)-g(\ba_2)|\over|\ba_1-\ba_2|}
\le|f|_\alpha|g|_\alpha,\label{HolProdBou}\en
as required.

Having established the two key properties of the H\"older norm, we can now prove
analyticity along the same lines as before: we introduce the generating function
\UB\zeta_\alpha(t)\equiv\sum_{s=1}^\infty|\nL\bxi^{(s)}|_\alpha\,t^s,\UE{holderGenfu}
where we use the notation $|\nL{\bf f}|_\alpha\equiv\max_{i,j}|\nL_if_j|_\alpha$.
Applying \rf{HolCijBou} and \rf{HolProdBou} to \rf{bigXis}, we deduce the inequality
\UB 6\Theta_\alpha\zeta^3_\alpha+12\Theta_\alpha\zeta^2_\alpha-\zeta_\alpha+T\ge0,\UE{HolderPoly}
where $T=t\,|\nL\vz|_\alpha$. Like in Section~\ref{sss:elementary}, this implies
analyticity of the Lagrangian map until a time inversely proportional
to the H\"older norm of the initial vorticity.

It is well-known that absolute summability of the Fourier series and H\"older
regularity define distinct classes of functions, although the difference is
subtle. For example, if a function is H\"older-continuous of exponent
$\alpha>1/2$, a theorem of Bernstein implies the absolute convergence
of the Fourier series \citep[see, e.g.,][Section~6.1]{Ka}. This result is sharp:
for $\alpha=1/2$, there exist functions that are H\"older-continuous without
having an absolutely summable Fourier series; a classical example is
the Hardy--Littlewood series $\sum_{n\ge1}n^{-1}\e^{\i Cn\ln n}\e^{\i nq}$ that
converges uniformly in the interval $[0,2\pi]$ to a H\"older-continuous function
of exponent 1/2 without being absolutely summable \citep[see, e.g.,][p.~197]{Zy}.

\subsection{Bounds on the radius of convergence of the time-Taylor series}
\label{ss:bounds}

After we have established the analyticity in time of the Lagrangian
trajectories under various assumptions on the initial conditions, it is
natural to ask: what is the radius of convergence of the time-Taylor
series that we have constructed?
This radius of convergence tells
us how far we can evaluate the solution without having to perform
analytic continuation. This is important when investigating
the issue of finite-time blow up for solutions of the incompressible
Euler equations. For example, \citet{Bra} investigated the blow up
in the Taylor--Green flow \citep{TG} and showed
by analysing the Taylor series for the enstrophy (the space integral of
the squared vorticity) that its radius of convergence is determined by
a pair of pure imaginary singularities at
the complex time locations $t_\star\approx\pm 2.18\i$.
Any putative real singularity has to be beyond the disk of convergence,
but in practice it is very difficult to analytically continue
outside the disk of convergence.

We have by now derived {\it lower
bounds} for the radius of convergence of the time-Taylor series of the
Lagrangian map. Can we obtain also {\it upper bounds}? For example,
can we ascertain that the Lagrangian map is not an entire function with infinite
radius of convergence? In general, we have not succeeded in doing this.
However, we have constructed analytically the full Lagrangian map and found
the radius of convergence to be finite for a steady-state solution to the Euler
equations for incompressible fluid. As any Beltrami flow, the well-known
ABC flows \citep[see][and references therein]{Dom}
are steady-state solution to the 3D Euler equations. In general, they have
a non-trivial Lagrangian structure, and the ordinary
differential equations (ODEs) for the Lagrangian trajectories are not analytically
integrable unless one of the coefficients $A$, $B$ or $C$ vanishes
(see Appendix A of \citet{Dom} and \citet{PaMa}). Here, we consider
the symmetric ``AB flow''
\UB\bv(x_1,x_2)=(-\sin x_1\sin x_2,\,-\cos x_1\cos x_2,\,\sqrt2\sin x_1\cos x_2).\UE{weirdAB}
It is easily checked that its vorticity is $\bom(x_1,x_2)=\sqrt2\,\bv(x_1,x_2)$
and that it becomes an ABC flow with $C=0$ and $A=B=1$ after a rotation
by $\pi/4$ about the $x_3$ axis, a stretching of the coordinates by a factor
$\sqrt2$ and a suitable permutation of the coordinates. In Section~\ref{sss:AB},
we calculate the radius of convergence of its Lagrangian map.

The next step (Section~\ref{sss:improvements}) is to improve the
lower bounds on the radius of convergence. While the numbers
obtained in Section~\ref{sss:elementary} can be improved by various methods,
particularly in the 2D case, as we shall see, a substantial gap remains between
the improved estimates and the actual values for the AB flow \rf{weirdAB}.

\subsubsection{AB flow: a flow with integrable Lagrangian trajectories}
\label{sss:AB}

The AB flow \rf{weirdAB} is independent of the vertical, $x_3$, coordinate.
As a consequence, the vertical component of the velocity behaves
just as a scalar advected by the horizontal components. The stream function of
the horizontal part of the flow is
\UB\psi(x_1,x_2)=\sin x_1\cos x_2,\UE{stream}
which, like the flow itself, is time-independent in Eulerian coordinates.
The Lagrangian map is as usual the solution to $\dot\bx=\bv$ with the initial
condition $\bx(0)=\ba$. It is easily checked that the stream function
is a Lagrangian invariant. Consequently, the third
component of the Lagrangian is just a linear function of time and cannot
induce loss of analyticity. Henceforth we consider only the horizontal
part of the fluid particle motion, governed by
\UB\dot x_1 =-\sin x_1\sin x_2,\qquad \dot x_2=-\cos x_1\cos x_2.\UE{horiz}

\begin{figure}
\centerline{\includegraphics[height=2in]{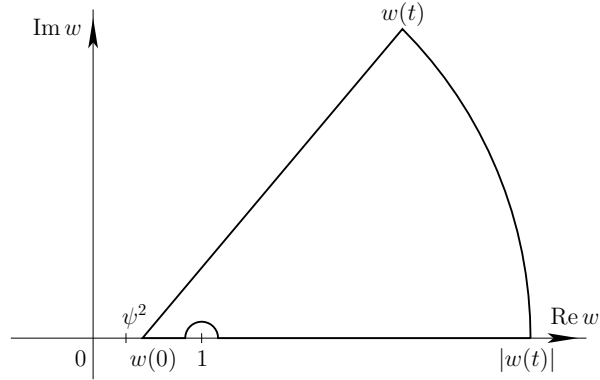}}
\caption{A sketch of the path of integration in the $\C$ plane
for the integral \rf{wInt} for computation of the minimal (in absolute value)
time of formation of singularity in the AB flow.}
\label{fig2}\end{figure}

To solve \rf{horiz} we make use of the invariance of the stream
function \rf{stream}. Introducing the new dependent variable
$w\equiv\sin^2x_1$, and assuming without loss of generality that
$\sin x_1\sin x_2\ge0$, we reduce \rf{horiz} to a separable ODE
\UB\dot{w}= -2\sqrt{w(1-w)(w-\psi^2)}.\UE{wEq}
Thus, the solution can be represented in terms of the elliptic integral:
\UB\int_{w(0)}^{w(t)}{\d w\over\sqrt{w(1-w)(w-\psi^2)}}=-2t.\UE{wInt}
For our purposes, we do not need much of the theory of elliptic
functions. The important point is that the only singularities
of the elliptic function $w(t)$ are double poles \citep[see, e.g.,][]{Akh}.
When $t$ approaches a time of singularity $t_\star$, the dependent variable $w$
tends to infinity and hence the original variable $x_1$ tends to complex
infinity. To calculate such singular times we should evaluate the definite integral
\rf{wInt} from the real value $w(0)$ (between $\psi^2$ and 1) to complex
infinity. Consider the contour shown in Fig.~\ref{fig2}, designed to avoid
inside it singularities of the integrand in \rf{wInt}. It is easily shown
that the integral over the arc from $w(t)$ to $|w(t)|$ tends to
zero as $w(t)$ tends to complex infinity. Letting the
radius of the semi-circle around $w=1$ tend to zero, we can replace the
integration from $w(0)$ to complex infinity by an integration over
the real axis from $w(0)$ to $+\infty$:
\UB2t_\star=\int_{w(0)}^1{\d w\over\sqrt{w(1-w)(w-\psi^2)}}
+\i\int_1^\infty{\d w\over\sqrt{w(w-1)(w-\psi^2)}}.\UE{timeInfty}

Next, we want to find the values $\psi$ and $w(0)$ that minimize
$|t_\star|$. By \rf{timeInfty}, the imaginary part of $t_\star$ does not depend
on $w(0)$. Hence the smallest module, for a given $\psi$, is obtained
for the smallest real part, namely, for $w(0)=1$. Finally, we have
to minimize over $\psi$. Inspection of the second, pure imaginary, integral
in \rf{timeInfty} shows that the minimum is achieved for $\psi=0$.
For these parameter values the integral becomes elementary:
\UB t_\star={\i\over2}\int_1^\infty{\d w\over w\sqrt{w-1}}={\i\pi\over2}\UE{mintime}
(it can be evaluated by substituting $w=1+\tan^2\theta$ for $0\le\theta\le\pi/2$).
Hence the radius of convergence of the time-Taylor series
for the AB flow is the finite number $\pi/2$.

\subsubsection{Improvements of bounds for the radius of convergence of the Taylor series}
\label{sss:improvements}

We shall now show that the lower bound for the radius of convergence
of the time-Taylor series for the Lagrangian map can be somewhat improved.
For this we shall use the Fourier series formulation and notation of Section~\ref{sss:elementary}. In that section, bounds on the Fourier coefficients of the
various time-Taylor coefficients were obtained from the recurrence relations
\rf{CauchyRecurrOpt}--\rf{sdiffRecurr} and then used to bound
the generating function $\zeta$ \rf{genfu}. Actually, the first few time-Taylor
coefficients can be calculated explicitly from \rf{CauchyRecurrOpt}--\rf{sdiffRecurr}
and then used to obtain tighter bounds.

For instance, for $s=1$ and $s=2$, we obtain
\UB\widehat{\bxi}^{(1)}_{\bf p}={\i{\bf p}\times\widehat{\bom}_{\bf p}\over|{\bf p}|^2}; \qquad\widehat{\bxi}^{(2)}_{\bf p}=-{\i{\bf p}\over2|{\bf p}|^2}
\sum_{0\ne{\bf r}\ne\bf p}{[{\bf r},{\bf p}-{\bf r},\widehat{\bom}_{\bf r}]
[{\bf r},{\bf p}-{\bf r},\widehat{\bom}_{\bf p-r}]\over
|{\bf r}|^2|{\bf p-r}|^2}.\UE{xi-one-and-two}
Here, the Fourier coefficients of the initial
vorticity are denoted, as before, $\widehat{\bom}_{\bf p}$.
Therefore, the bounds \rf{FCauchyTFbound}--\rf{sdiffTFbound}
for the transverse and longitudinal components, and \rf{xi-one-and-two}
yield the following bound for the generating function:
\be\zeta&\le&t\sum_{\bf p}|{\bf p}|\,|\widehat{\bxi}^{(1)}_{\bf p}|+
t^2\sum_{\bf p}|{\bf p}|\,|\widehat{\bxi}^{(2)}_{\bf p}|
+\sum_{s=3}^\infty\sum_{\bf p}\left(|{\bf p}\times\widehat{\bxi}^{(s)}_{\bf p}|
+|{\bf p}\cdot\widehat{\bxi}^{(s)}_{\bf p}|\right)t^s\nonumber\\
&\le&T+{t^2\over2}\sum_{{\bf p}\ne0}\left|\sum_{0\ne{\bf r}\ne\bf p}{[{\bf r},{\bf p}-{\bf r},
\widehat{\bom}_{\bf r}]\,[{\bf r},{\bf p}-{\bf r},\widehat{\bom}_{\bf p-r}]
\over|{\bf r}|^2|{\bf p-r}|^2}\right|\nonumber\\
&&+\,{1\over\sqrt2}\sum_{s=3}^\infty\sum_{0<m<s,\,\bf r}|{\bf r}|\,|{\bf p}-{\bf r}|
|\widehat{\bxi}^{(m)}_{\bf r}|\,|\widehat{\bxi}^{(s-m)}_{{\bf p}-{\bf r}}|\,
t^s\nonumber\\
&&+{1\over6}\sum_{s=3}^\infty\,\sum_{\subs{n_1+n_2+n_3=s}{{\bf r}_1+{\bf r}_2+{\bf r}_3=\bf p}}
|{\bf r}_1|\,|{\bf r}_2|\,|{\bf r}_3|\,|\widehat{\bxi}^{(n_1)}_{{\bf r}_1}|\,
|\widehat{\bxi}^{(n_2)}_{{\bf r}_2}|\,
|\widehat{\bxi}^{(n_3)}_{{\bf r}_3}|\,t^s\nonumber\\
&=&T+{t^2\over2}\sum_{{\bf p}\ne0}\left|\sum_{0\ne{\bf r}\ne\bf p}{[{\bf r},{\bf p}-{\bf r},\widehat{\bom}_{\bf r}]
[{\bf r},{\bf p}-{\bf r},\widehat{\bom}_{\bf p-r}]\over|{\bf r}|^2|{\bf p-r}|^2}\right|
+{\zeta^2-T^2\over\sqrt2}+{\zeta^3\over6}\label{better}\\
&\le&T+T^2\left({1\over2}-{1\over\sqrt2}\right)+{\zeta^2\over\sqrt2}+{\zeta^3\over6}.
\label{lessbetter}\en
Here, $T\equiv t\Gamma = t\sum_{\bf p}|\widehat{\bom}_{\bf p}|$, as before.
Observe that \rf{lessbetter} differs from \rf{final} only by a new
quadratic term in the dimensionless time $T$. It immediately follows that
convergence of the series
\rf{xiSeries} is guaranteed for $t<\widetilde t_c\equiv\widetilde T_c/\Gamma$,
where the critical dimensionless time $\widetilde T_c$ satisfies the equation
$\widetilde T_c+(1/2-1/\sqrt{2})\widetilde T^2_c=T_c$ implying
\UB\widetilde T_c=1+\sqrt{2}-\sqrt{13/3}\approx0.3325\,.\UE{newTcrit}
Comparison with \rf{Tcrit} show a modest improvement of about $7\%$.

So far we were assuming an arbitrary 3D initial condition (within the class
of vorticities with summable Fourier coefficients). A more significant improvement can
be achieved by constructing bounds for a specific flow in two dimensions.
For instance, for the AB flow \rf{stream}, for which we know the exact radius
of convergence, the Taylor coefficients of order one and two are easily shown to be
\UB\bxi^{(1)}=(-\sin a_1\sin a_2,-\cos a_1\cos a_2),\qquad \bxi^{(2)}={1\over4}\,(\sin2a_1,-\sin2a_2).\UE{ABxionetwo}
Consequently, for the AB flow,
\UB\Gamma=\sum_{\bf p}|{\bf p}|\,|\widehat{\bxi}^{(1)}_{\bf p}|=2,
\qquad\sum_{\bf p}|{\bf p}|\,|\widehat{\bxi}^{(2)}_{\bf p}|=1,\UE{ABxiOneTwo}
and thus the analogue of \rf{better} takes the form
\UB T+T^2\left({1\over4}-{1\over\sqrt2}\right)-\zeta+{\zeta^2\over\sqrt2}\ge0\UE{ABbou}
(in two dimensions the cubic term is absent in \rf{better}).
In the dimensionless time $T$, analyticity holds now until
\UB T_{\rm AB}={2-\sqrt[4]{2}\over\sqrt{8}-1}\approx0.4434,\UE{TcritAB}
an improvement of about $43\%$ over the value given by \rf{Tcrit}.
The actual value of the radius of convergence for the
AB flow in dimensionless form is $\pi$, which is still about a factor 7 larger
than our improved bound. There is definitely room for further improvement.

\section{Ideal compressible fluid flow in an Einstein--de Sitter universe}
\label{s:cosmo}

The simplest model for the dynamics of compressible ideal fluid
in the potential case, in any dimension, is the Burgers equation:
\UB\partial_t\bv+(\bv\cdot\nabla)\bv=0;\qquad\bv =-\nabla\varphi.\UE{burgers}
Its Lagrangian formulation is simply $\ddot\bx=0$, which
integrates to $\bx(\ba,t)=\ba+t\vz(\ba)$. This Lagrangian map is
obviously analytic, and even entire, in the time variable. Singularities
appear nevertheless, but only when reverting to Eulerian coordinates,
because the Lagrangian map ceases to be invertible when its
Jacobian $J\equiv\det(I+t\nL\vz)$ vanishes for the first time (here
$I$ denotes the identity matrix).

There are models of compressible dynamics with far less trivial
Lagrangian maps that can be viewed as extensions of the ideal potential Burgers
dynamics. Of particular interest are the Euler--Poisson equations, briefly
presented in Section~\ref{sss:cosmoEuler}, which arise in cosmology and govern
the essentially collisionless and thus ideal motion of dark matter.
Actually, our interest in analyticity of Lagrangian trajectories in hydrodynamics was
triggered by work on the Lagrangian perturbation theory in cosmology
\citep[][see also \citet{Bu89}]{Mou,Bo92,Bo95,Bu92}.

\subsection{The Euler--Poisson equations}\label{ss:EP}

\subsubsection{Basic equations in Eulerian coordinates}\label{sss:cosmoEuler}

The cosmological Euler--Poisson system of equations takes its shortest
form in the so-called comoving coordinates $\bx$ (moving with the
Hubble expansion). We denote by $\bv$ the peculiar velocities (with
the Hubble expansion subtracted), by $\rho$ the normalised matter
density, by $\phig$ the gravitational potential and by $\tau$ a
suitable time variable (see below). Since, in this paper, we want to
concentrate on fluid dynamical aspects, we shall limit ourselves to an
Einstein--de~Sitter (EdS) model of a flat matter-dominated universe,
ignoring the cosmological term responsible for the accelerated
expansion of the Universe. For the EdS model the Euler--Poisson equations take
the following form \citep[cf.][and references therein]{MNRAS}:
\be\partial_\tau\bv+(\bv\cdot\nabla)\bv&=&-{3\over2\tau}(\bv+\nabla\phig),\label{comoveuler}\\
\partial_\tau\rho+\nabla\cdot(\rho\bv)&=&0,\label{comovcontinuity}\\
\nabla^2_{\rm E}\phig&=&{\rho-1\over\tau}\label{comovpoisson}\en
(here $\nabla^2_{\rm E}$ is the Eulerian Laplacian).
Eq.~\rf{comoveuler} is the momentum equation,
\rf{comovcontinuity} expresses mass conservation, \rf{comovpoisson} is
the Poisson equation for the gravitational potential. One advantage
of the variables here chosen is that such cosmological equations
significantly resemble standard fluid mechanical equations.

We comment now on the time variable $\tau$, which is not the
cosmic time $t$, but is proportional to $t^{2/3}$. Let us linearise
the Euler--Poisson equations about the steady state $\bv=0,\, \rho=1$. Introducing the density perturbation
$\delta\equiv \rho -1$, we obtain from \rf{comoveuler}--\rf{comovpoisson}
\UB\partial^2_{\tau\tau}\delta=-\frac{3}{2\tau}\left(\partial_\tau\delta-\frac{\delta}{\tau}\right).\UE{linearperturb}
It is immediately seen that \rf{linearperturb} has two power-law solutions,
$\tau$ and $\tau^{-3/2}$. In the former, density perturbations
grow linearly with $\tau$, hence the name ``linear growth time'' for the
variable $\tau$. (The standard notation in cosmology for the linear growth
time is $D$ or $D(t)$.) The latter mode is
decaying at large times but blows up as $\tau\to 0$, thus invalidating
linearisation. In the full cosmological setting $\tau$ may not go all the way
to zero, because the matter-dominated
description used for deriving the Euler--Poisson equations is valid only after
the decoupling of matter and radiation. In the cosmic time variable $t$,
starting with the Big Bang,
this happens around $t=380,000$ years, a kind of boundary layer when
compared to the present age of the Universe $t\approx 13.8$ billion years.
We shall ignore such effects and let $\tau$ become arbitrarily small. Since
$\tau$ appears in the denominator of the right-hand sides of \rf{comoveuler}
and \rf{comovpoisson}, if we demand that the solutions remain well-behaved
as $\tau\to0$, the initial conditions must satisfy two
{\it slaving conditions}, $\rho\to1$ and $\bv\to-\nabla\phig^{\rm(init)}$
for $\tau\to0$. Eq.~\rf{comoveuler} then implies that the velocity remains
potential in Eulerian (but not Lagrangian!) coordinates at any $\tau>0$.

\subsubsection{Basic equations in Lagrangian variables}\label{sss:cosmoLag}

Fluid dynamicists, since the time of Lagrange, denote what was later
called Lagrangian coordinates by the letter $\ba$. We shall follow this
tradition. Cosmologist typically denote it by $\bm q$, to avoid the letter $a$,
standing for the expansion scale factor. We, however, avoid the sub- or
superscript 0 when referring to initial conditions --- in cosmology zero
refers to zero redshift, that is the present epoch. Instead, as in
Section~\ref{s:ideal}, we use the superscript (init).

To derive the Lagrangian form of the Euler--Poisson equations from
\rf{comoveuler}--\rf{comovpoisson}, we shall use a mixture of the
methods already employed for the incompressible case, and of a technique
widely used in cosmological fluid dynamics \citep[see, e.g.,][]{EB97}.

As in the incompressible case \rf{EulerMom}, because the velocity
is potential, the right-hand side of \rf{comoveuler} is again a Eulerian
gradient. Thus we can apply identically the Cauchy--Weber derivation
of Section~\ref{ss:cauchyinv}. The only difference is that now the initial
vorticity vanishes. Thus, \rf{cauchy} reduces to
\UB\sum_{k=1}^3\nL\dot{x}_k\times\nL x_k=0.\UE{cosmocauchy}
A related cubic equation, expressing the vanishing of
the current Eulerian vorticity is found in \citet[][cf.~eq.~(2.24)]{RB12}.

The second equation, replacing the incompressibility condition of
the hydrodynamic problem, stems from mass conservation
expressed as $\rho J=1$, where $J=\det(\nL x)$ is
the Jacobian of the Lagrangian map. We note that the left-hand side
of \rf{comoveuler} is the acceleration $\ddot\bx$. Taking the Eulerian
divergence of \rf{comoveuler}, using \rf{comovpoisson} and $\bv=\dot\bx$,
we obtain
\UB\nabla\cdot(\ddot\bx+{3\over2\tau}\dot\bx)=
-{3\over2\tau^2}\nabla^2_{\rm E}\phig=-{3\over2\tau^2}(\rho-1).\UE{precosmodiv}
Next, we have to reexpress the Eulerian divergence in terms
of Lagrangian space derivatives. For this, we need
the inverse of the Jacobian matrix. We use the following identity, true
for an arbitrary smooth 3D map:
\UB\nabla_ia_j={\varepsilon_{jmn}\,\varepsilon_{ikl}\over2J}\,\nL_mx_k\,\nL_nx_l.\UE{dqldxi}
(This expression is more convenient than the one involving cofactor matrices.)
Eq.~\rf{precosmodiv} becomes then
\UB\varepsilon_{jmn}\,\varepsilon_{ikl}\,\nL_mx_k\,\nL_nx_l\,
\left(\partial^2_{\tau\tau}+{3\over2\tau}\,\partial_\tau\right)\nL_jx_i
={3\over\tau^2}(J-1).\UE{cosmosdiff}
Eq.~\rf{cosmocauchy} and \rf{cosmosdiff} constitute a closed system of
Lagrangian equations for fluid particle motion. The equations, particularly
\rf{cosmosdiff}, may have solutions possibly singular at $\tau=0$ due
to the presence of negative powers of $\tau$. Actually, with the slaving
conditions (or, equivalently, the absence of the decaying mode), it is
easily checked that solutions exist in the form of power series
that are non-singular at $\tau=0$. Such expansions in powers of the linear
growth time $\tau$ are widely used in the cosmological literature on the Lagrangian perturbation expansion for an EdS universe \citep[see, e.g.,][]{Mat08}. The formal power series is
constructed in Section~\ref{ss:cosmoTaylor}. Convergence and analyticity
are then established in Section~\ref{ss:cosmoHolder}.

We also observe that, as in the incompressible case, the system
\rf{cosmocauchy} and \rf{cosmosdiff} is invariant under a change to an
arbitrarily accelerated frame. In cosmology this is called
``extended Galilean invariance'' \citep[see][and references therein]{Ber13}.

\subsection{The formal Taylor expansion}\label{ss:cosmoTaylor}

We now seek a solution to \rf{cosmocauchy}--\rf{cosmosdiff} in the form
of a power series in $\tau$ for the displacement $\bxi\equiv\bx-\ba$:
\UB\bxi(\ba,t)=\sum_{s=1}^\infty\bxi^{(s)}(\ba)\tau^s.\UE{cosmoSeries}
Upon substitution of this series into \rf{cosmosdiff} we
obtain the following relations among the Taylor coefficients
(summation over repeated indices is now assumed):
\be(2s+3)(s-1)\nL\cdot\bxi^{(s)}&=&\!\!\!\!\sum_{n_1+n_2 =s}\!\!\left(2n_1^2+n_1-{3\over2}\right)\left(\nL_j\xi_i^{(n_2)}\,\nL_i\xi_j^{(n_1)}
-\nL_i\xi_i^{(n_2)}\,\nL_j\xi_j^{(n_1)}\right)\nonumber\\
&-&\!\!\!\sum_{n_1+n_2+n_3=s}\left(n^2_3+{n_3-1\over2}\right)\varepsilon_{jmn}\,\varepsilon_{ikl}\,\nL_m\xi^{(n_1)}_k\,\nL_n\xi^{(n_2)}_l\,\nL_j\xi^{(n_3)}_i\,.\nonumber\\
\label{preCosmoSdiffRecurr}\en
This equation is trivially satisfied for $s=1$. For $s>1$, we symmetrise
both sums in the right-hand side by performing all possible permutations
of indices to obtain
\be\nL\cdot\bxi^{(s)}\!
&=&\!\!\sum_{0<n<s}\!\!{n^2+(s-n)^2+(s-3)/2\over s^2+(s-3)/2}\!\!\sum_{1\le i<j\le3}\!
\left(\nL_j\xi_i^{(s-n)}\,\nL_i\xi_j^{(n)}\!-\nL_i\xi_i^{(s-n)}\,\nL_j\xi_j^{(n)}\right)\nonumber\\
&&\!\!\!-{1\over6}\sum_{n_1+n_2+n_3=s}\!\!\!{n^2_1+n^2_2+n^2_3+(s-3)/2\over s^2+(s-3)/2}\,\varepsilon_{jmn}\,\varepsilon_{ikl}\,\nL_m\xi^{(n_1)}_k\,\nL_n\xi^{(n_2)}_l\,\nL_j\xi^{(n_3)}_i\,.\nonumber\\
\label{cosmoSdiffRecurr}\en
The other relations, stemming from \rf{cosmocauchy}, are obtained
for $s\ge 1$ just as in Section~\ref{ss:taylor}:
\UB\nL\times\bxi^{(s)}={1\over2}\sum_{\subs{1\le k\le3,}{0<n<s}}{s-2n\over s}\,
\nL\xi^{(n)}_k\times\nL\xi^{(s-n)}_k\,.\UE{cosmoCauchyRecurr}

Consequently, for $s\ge1$ the Taylor coefficients $\bxi^{(s)}$ can be
determined by performing a Helmholtz--Hodge decomposition \rf{xis}
(now, ${\bf R}^{(1)}$ and $R^{(2)}$ denote the right-hand sides
of \rf{cosmoCauchyRecurr} and \rf{cosmoSdiffRecurr}, respectively).

Observe that the recurrence relations \rf{cosmoSdiffRecurr}--\rf{cosmoCauchyRecurr}
for this compressible case turn out to be identical in structure to
\rf{CauchyRecurrOpt}--\rf{sdiffRecurr} for the incompressible one. This is
not surprising given that, in Lagrangian coordinates, the fluid particle
displacements have rotational and gradient components in both cases.

In the cosmological literature similar expansions are found, which are
not always time-Taylor expansions, but expansions in powers of the magnitude
of the displacement. They have been carried out up to the fourth order \citep{RB12}.
To the best of our knowledge, explicit all-order recurrence relations
such as \rf{cosmoSdiffRecurr}--\rf{cosmoCauchyRecurr} have never been presented
in the literature. They are essential if we want to investigate the
convergence of the expansion and the analyticity of Lagrangian trajectories.

We have checked that our recurrence relations are consistent with the known
results up to the fourth order. In particular, the first term $\bxi^{(1)}=\nL\phig^{\rm(init)}$
yields the Burgers--Zeldovich approximation, in which fluid particles
keep their initial velocities, often used to initialise numerical simulations.
It is important to observe that, in the one-dimensional case when
the initial velocity depends only on a single coordinate, the Taylor
expansion terminates exactly with the term linear in $\tau$.
In other words, in one dimension the Zeldovich approximation is exact
\citep{Bu92}.

\subsection{Analyticity in time of fluid-particle trajectories}\label{ss:cosmoHolder}

The Lagrangian perturbation theory in cosmology has been, since its introduction about
twenty years ago, a formal expansion in powers of a parameter controlling
the amplitude of particle displacements. For the case of an EdS
universe, the linear growth time $\tau$ is an appropriate
control parameter and the Lagrangian map becomes a formal
Taylor series in powers of $\tau$. (Analyticity issues for non-EdS models having
for example a cosmological term, are beyond the scope of the present paper.)
With the tools developed in
Section~\ref{s:ideal} it is now easy to see that with suitable
conditions on the initial velocity field (or, equivalently, on the initial
gravitational potential) this formal expansion is actually convergent
for small enough $\tau$ and defines an analytic function of $\tau$.
Indeed, given the structure of the recurrence relations
\rf{cosmoSdiffRecurr}--\rf{cosmoCauchyRecurr}, whenever the initial velocity
gradient is in any of the function spaces considered in Sections~\ref{ss:L1},
\ref{ss:analytic} and \ref{ss:holder}, analyticity will be guaranteed
for at least a finite time. To just show analyticity, it suffices to
observe that all
the coefficients appearing in \rf{cosmoSdiffRecurr}--\rf{cosmoCauchyRecurr}
are uniformly bounded by unity. Improved bounds on the radius of convergence
of the time-Taylor series can be obtained for the case of an initial
velocity gradient with a summable Fourier series just as in
Section~\ref{sss:elementary}: all the estimates derived
in the incompressible case still hold, without any alterations,
in the compressible Euler--Poisson case.
Further enhancements, such as discussed in Section~\ref{sss:improvements}
can also be implemented (and actually yield a larger radius of convergence
of the Taylor series). For the case of an initial velocity analytic in the
space variable, the proof of Section~\ref{ss:analytic} that spatio-temporal
analyticity persists for at least a finite time holds identically.

Let us illustrate this for the simple case when the initial velocity gradient
$\nabla^2\phig^{\rm(init)}$ is H\"older-continuous of exponent $\alpha$.
As in Section~\ref{ss:holder}, we introduce the generating function
\UB\zeta_\alpha(\tau)\equiv\sum_{s=1}^\infty|\nL\bxi^{(s)}|_\alpha\,\tau^s.\UE{cosmoGenfu}
Applying \rf{xis}, \rf{cosmoSdiffRecurr} and \rf{cosmoCauchyRecurr}, we find
\UB6\Theta_\alpha\zeta_\alpha^3+12\Theta_\alpha\zeta_\alpha^2-\zeta_\alpha+\tilde\tau\ge0,\UE{cosmoPoly}
where $\tilde\tau=\tau\max_{ij}|\nL_i\nL_j\phig^{\rm(init)}|_\alpha\,$. Proceeding
as in \citet{FriZh}, we find that analyticity holds until $\tau_c={\cal T}/\max_{ij}|\nL_i\nL_j\phig^{\rm(init)}|_\alpha$, where for small $\alpha$
\UB{\cal T}=(48\Theta_\alpha)^{-1}+{\rm O}(\Theta^{-2}_\alpha)\sim\alpha\,.\UE{qc}
Thus, the smaller the H\"older exponent $\alpha$ of the initial velocity
gradient, the smaller is the guaranteed interval of time analyticity.
If however the initial velocity gradient, or the initial density fluctuations
are not sufficiently smooth, the series may well fail to converge.
For example, \citet{SS96} investigated the case of an initial
``top-hat underdensity'' that is initially discontinuous, and found that
low-order perturbation worked better than higher-orders, which they regarded
as a possible evidence for ``semiconvergence'' of the perturbation series.

\section{Concluding remarks}\label{s:conclusion}

In this paper we have considered two instances of ideal three-dimensional
flow, incompressible flow (Section~\ref{s:ideal}) and a gravity-driven
compressible flow of cosmological relevance (Section~\ref{s:cosmo}). In both
cases, we have shown that the motion of fluid particles is time-analytic,
for at least a finite time, when the initial vorticity or velocity gradient,
respectively, are slightly better than just continuous in the space variables
(for example, H\"older-continuous).

Ours are solutions to the Euler or the Euler--Poisson problems with a
Cauchy-type formulation in Lagrangian coordinates. But do they yield
solutions to the equations in the original formulations in Eulerian
coordinates? The latter are obtained from the former ones by using the
global inverse Lagrangian map, whose existence is therefore required
for such a construction. For the case of an incompressible flow, the
existence of a global inverse Lagrangian map is easily proved using
the conservation of volumes. In the compressible Euler--Poisson case,
analyticity until some time $\tau_\star$ does not imply existence of such
a global inverse Lagrangian map. Indeed, in the particular case of one
dimension, for which the Burgers--Zeldovich approximation is exact
(cf.~Section~\ref{ss:cosmoTaylor}), the Lagrangian map is a linear
function of $\tau$ that is trivially analytic for all $\tau$, but the
invertibility holds only until the time of the first crossing of
particle paths (called ``shell crossing'' in cosmology).

What happens to fluid particles beyond the real positive time, say,
$t_1$, until which our results guarantee analyticity of their
trajectories? If the invertibility requirement for the Eulerian
formulation does hold till $t_1$ and the Eulerian solution is still
in the appropriate function space (for example, has a vorticity whose
Fourier series is absolutely summable), then we can again obtain
analytical trajectories till time $t_2>t_1$ by constructing a new
time-Taylor series, starting at time $t_1$ (in the cosmological
setting, minor modifications are required, since \rf{cosmosdiff} is
not autonomous in time). This procedure can be iterated any number of
times until one of the conditions stated above is violated and the
Eulerian solution becomes therefore rougher than initially.

Does this happen and if so when? For the case of incompressible flow, this is
actually the issue of whether a blow up of solutions can occur after a finite
time, a major mathematical question. By contrast, for compressible flow
in cosmology there is hardly any doubt that such finite-time blow up does occur,
essentially because fast particles can catch up with slower ones. In this
case, however, a physical question is open: is there an interval of time after
the matter-radiation decoupling, during which no particle crossing occurs? This
depends very much on what are the small-scale properties of the density and
velocity fluctuations at decoupling that are inherited from primordial
cosmology, questions that will be addressed elsewhere.

Finally, we want to stress again how remarkably similar are the Lagrangian
dynamics of incompressible and compressible ideal flow when tackled in
the spirit of Cauchy's 1815 work.

We are grateful to C.~Bardos, F.~Bernardeau, F.~Bouchet, T.~Buchert, S.~Colombi,
J.-P.~Kahane, G.~Lebeau, H.K.~Moffatt, A.~Shnirelman, L.~Sz\'ekelyhidi
and B.~Villone for fruitful discussions. VZ was supported
by the grant 11-05-00167-a from the Russian foundation for basic research.
His visits to the Observatoire de la C\^ote d'Azur (France) were
supported by the French Ministry of Higher Education and Research.

\appendix
\section{}\label{appA}

Each step in the derivation of \rf{sdiffTF} is elementary, but their sequence
is not necessarily intuitive; this has prompted us to present here
the derivation in detail.

Upon substituting the Fourier series \rf{TaylorFourier}
into \rf{sdiffRecurr} and collecting the terms involving the exponent
$\e^{\i{\bf p}\cdot\ba}$ we find
\be{\bf p}\cdot\widehat{\bxi}^{(s)}_{\bf p}\!&=&-\i\sum_{0<m<s,\ \bf r}\ \sum_{1\le i<j\le3}
\varepsilon_{ijk}[{\bf r},{\bf p}-{\bf r},{\bf e}_k]
\left(\widehat{\bxi}^{(m)}_{\bf r}\right)_{\!i}
\left(\widehat{\bxi}^{(s-m)}_{{\bf p}-{\bf r}}\right)_{\!j}\nonumber\\
&&+\sum_{\subs{n_1+n_2+n_3=s}{{\bf r}_1+{\bf r}_2+{\bf r}_3=\bf p}}
[{\bf r}_1,{\bf r}_2,{\bf r}_3]\,
\left(\widehat{\bxi}^{(n_1)}_{{\bf r}_1}\right)_{\!1}
\left(\widehat{\bxi}^{(n_2)}_{{\bf r}_2}\right)_{\!2}
\left(\widehat{\bxi}^{(n_3)}_{{\bf r}_3}\right)_{\!3},\label{presdiffTF}\en
where ${\bf e}_k$ denotes the $k$-th unit vector of the cartesian coordinate
system, $(\widehat\bxi_{\bf r})_j$ the $j$-th component of the vector
$\widehat\bxi_{\bf r}$, and $k\equiv6-i-j$. The right-hand side
of \rf{presdiffTF} consists of two sums.
We replace the first one by the arithmetic mean of this sum and its copy, in which
the pairs of indices $(m,{\bf r})\leftrightarrow(s-m,{\bf p-r})$ are swapped.
Similarly, we replace the second sum by the arithmetic mean of this sum and
its five copies with permuted pairs of indices $(n_i,{\bf r}_i)$. The replacements
involve a group of two terms quadratic in $\widehat{\bxi}^{(n_i)}_{{\bf r}_i}$
and a group of six terms cubic in $\widehat{\bxi}^{(n_i)}_{{\bf r}_i}$; each
group readily amalgamates into a triple product, yielding
\be{\bf p}\cdot\widehat{\bxi}^{(s)}_{\bf p}
&=&-{\i\over2}\sum_{0<m<s,\ \bf r}\ \sum_{k=1}^3\,
[{\bf r},{\bf p}-{\bf r},{\bf e}_k]\,\left[\widehat{\bxi}^{(m)}_{\bf r},
\widehat{\bxi}^{(s-m)}_{{\bf p}-{\bf r}},{\bf e}_k\right]\nonumber\\
&&+{1\over6}\sum_{\subs{n_1+n_2+n_3=s}{{\bf r}_1+{\bf r}_2+{\bf r}_3=\bf p}}
[{\bf r}_1,{\bf r}_2,{\bf r}_3]\,\left[\widehat{\bxi}^{(n_1)}_{{\bf r}_1},
\widehat{\bxi}^{(n_2)}_{{\bf r}_2},\widehat{\bxi}^{(n_3)}_{{\bf r}_3}\right].\label{almsdiffTF}\en
Finally, we note that summation in $k$ in the first sum in \rf{almsdiffTF}
amounts to calculation of the scalar product; this gives rise to \rf{sdiffTF}.

\begin{thebibliography}{99}
\bibitem[Akhiezer(1990)]{Akh}
{\sc Akhiezer, N.I.} 1990 {\it Elements of the theory of elliptic functions}.
AMS Translations of Mathematical Monographs vol.~79. AMS, Rhode Island.

\bibitem[Arfken \& Weber(2005)]{AWe}
{\sc Arfken, G.B.~\& Weber, H.J.} 2005 {\it Mathematical method
for physicists}. Elsevier, Amsterdam.

\bibitem[Benachour(1976)]{Be76}
{\sc Benachour, S.} 1976 Analyticit\'e des solutions p\'eriodiques de l'\'equation
d'Euler en trois dimensions. {\it C.~R.~Acad.~Sci.~Paris} S\'erie A {\bf 283}, 107--110.

\bibitem[Bernardeau(2013)]{Ber13}
{\sc Bernardeau, F.} 2013 {\it The evolution of the large-scale structure
of the universe: beyond the linear regime.} Lectures given to the Les Houches
Summer School Post-Planck Cosmology, 8 July -- 2 August 2013
\url{http://arxiv.org/abs/1311.2724}

\bibitem[Bouchet {\it et al.}(1992)]{Bo92}
{\sc Bouchet, F.R., Juszkiewicz, R., Colombi, S.~\& Pellat, R.} 1992
Weakly nonlinear gravitational instability for arbitrary $\Omega$.
{\it Astrophys.~J.} {\bf 394}, L5--L8
\url{http://adsabs.harvard.edu/abs/1992ApJ...394L...5B}.

\bibitem[Bouchet {\it et al.}(1995)]{Bo95}
{\sc Bouchet, F.R., Colombi, S., Hivon, E.~\& Juszkiewicz, R.} 1995 Perturbative
Lagrangian approach to gravitational instability. {\it Astron.~Astrophys.}
{\bf 296}, 575--608 \url{http://adsabs.harvard.edu/abs/1995A%26A...296..575B}.

\bibitem[Brachet {\it et al.}(1983)]{Bra}
{\sc Brachet, M.E., Meiron, D.I., Orszag, S.A., Nickel, B.G., Morf, R.H.~\&
Frisch, U.} 1983 Small-scale structure of the Taylor-Green vortex.
{\it J.~Fluid Mech.} {\bf 130}, 411--452.

\bibitem[Brenier {\it et al.}(2003)]{MNRAS}
{\sc Brenier, Y., Frisch, U., H\'enon, M., Loeper, G., Matarrese, S., Mohayaee, R.~\&
Sobolevskii, A.} 2003 Reconstruction of the early Universe as a convex optimization
problem. {\it Mon.~Not.~R.~Astron.~Soc.} {\bf 346}, 501--524. arXiv:astro-ph/0304214.

\bibitem[Buchert(1989)]{Bu89}
{\sc Buchert, T.} 1989 A class of solutions in Newtonian cosmology and the
pancake theory {\it Astron. Astrophys.} {\bf 22}, 9--24.

\bibitem[Buchert(1992)]{Bu92}
{\sc Buchert, T.} 1992 Lagrangian theory of gravitational instability
of Friedman--Lemaitre cosmologies and the ``Zel'dovich approximation''
{\it Mon.~Not.~R.~Astron.~Soc.} {\bf 254}, 729--737
\url{http://adsabs.harvard.edu/abs/1992MNRAS.254..729B}.

\bibitem[Buckmaster {\it et al.}(2013)]{BuLeSz}
{\sc Buckmaster, T., De Lellis, C.~\& Sz\'ekelyhidi Jr, L.} 2013 Transporting
microstructure and dissipative Euler flows. arXiv:1302.2815 [math.AP]

\bibitem[Cauchy(1815)]{Cau27}
{\sc Cauchy A.L.} 1815 Sur l'\'etat du fluide \`a une \'epoque quelconque du mouvement.
M\'emoires extraits des recueils de l'Acad\'emie des sciences de l'Institut
de France, Th\'eorie de la propagation des ondes \`a la surface d'un fluide
pesant d'une profondeur ind\'efinie ({\it Extraits des M\'emoires pr\'esent\'es
par divers savants \`a l'Acad\'emie royale des Sciences de l'Institut de France
et imprim\'es par son ordre). Sciences math\'ematiques et physiques}. Tome I, 1827
Seconde Partie, pp.~33--73. \url{http://gallica.bnf.fr/ark:/12148/bpt6k90181x.r=Oeuvres+completes+d%27Augustin+Cauchy.langFR}

\bibitem[Chemin(1992)]{Ch92}
{\sc Chemin J.-Y.} 1992 R\'egularit\'e des trajectoires des particules
d'un fluide incompressible remplissant l'espace. {\it J.~Math\'ematiques Pures
et Appliqu\'ees} {\bf 71}, 407--417.

\bibitem[Constantin, Vicol \& Wu(2014)]{CVW}
{\sc Constantin, P., Vicol, V.~\& Wu, J.} 2014 Analyticity of Lagrangian
trajectories for well-posed inviscid incompressible fluid models.
Submitted. arXiv:1403.5749 [math.AP]

\bibitem[Dombre {\it et al.}(1986)]{Dom}
{\sc Dombre, T., Frisch, U., Greene, J.M., H\'enon, M., Mehr, A.~\& Soward, A.M.}
1986 Chaotic streamlines and Lagrangian turbulence: the ABC-flows.
{\it J.~Fluid Mech.} {\bf 167}, 353--391.

\bibitem[Ehlers \& Buchert(1997)]{EB97}
{\sc Ehlers, J.~\& Buchert, T.} 1997 Newtonian cosmology in Lagrangian
formulation: foundations and perturbation theory. {\it Gen. Rel. Grav.}
{\bf 29}, 733--764.

\bibitem[Frisch \& Villone(2014)]{FrVi}
{\sc Frisch, U.~\& Villone, B.} 2014 Cauchy's almost forgotten Lagrangian
formulation of the Euler equation for 3D incompressible flow. {\it European
Physical J.}~H, submitted. arXiv:1402.4957 [math.HO].

\bibitem[Frisch \& Zheligovsky(2014)]{FriZh}
{\sc Frisch, U.~\& Zheligovsky, V.} 2014 A very smooth ride in a rough sea.
{\it Comm.~Math. Physics} {\bf 326}, 499--505. arXiv:1212.4333 [math.AP].

\bibitem[Gilbarg \& Trudinger(1998)]{GT98}
{\sc Gilbarg, D.~\& Trudinger, N.S.} 1998 {\it Elliptic partial differential
equations of second order}. Springer, Berlin.

\bibitem[Hankel(1861)]{Han}
{\sc Hankel, H.} 1861 {\it Zur allgemeinen Theorie der Bewegung
der Fl\"ussigkeiten}, Preisschrift der philosophischen Facult\"at der Georgia
Augusta, G\"ottingen. \url{http://babel.hathitrust.org/cgi/pt?id=mdp.39015035826760;view=1up;seq=5}

\bibitem[H\"older(1933)]{Ho33}
{\sc H\"older, E.} 1933 \"Uber die unbeschr\"ankte Fortsetzbarkeit einer
stetigen ebenen Bewegung in einer unbegrenzten inkompressiblen Fl\"ussigkeit.
{\it Mathematische Zeitschrift\/} {\bf 37}, 727--738.

\bibitem[Isett(2014)]{Is}
{\sc Isett, Ph.} 2014 Regularity in time along the coarse scale flow
for the incompressible Euler equations. arXiv:1307.0565 [math.AP].

\bibitem[Katznelson(2004)]{Ka}
{\sc Katznelson, Y.} 2004 {\it An introduction to harmonic analysis}.
3rd ed., Cambridge Univ.~Press.

\bibitem[Korn(1907)]{Ko}
{\sc Korn, A.} 1907 Sur les \'equations de l'\'elasticit\'e. {\it Annales
scientifiques de l'\'E.N.S.} $3^{\rm e}$ s\'erie, {\bf 24}, 9--75.

\bibitem[Kukavica \& Vicol(2009)]{KV09}
{\sc Kukavica, I.~\& Vicol, V.} 2009 On the radius of analyticity of solutions
to the three-dimensional Euler equations. {\it Proc.~Am.~Math.~Soc.}
{\bf 137}, 669--677.

\bibitem[Kukavica \& Vicol(2011)]{KV11}
{\sc Kukavica, I.~\& Vicol, V.} 2011 On the analyticity and Gevrey-class
regularity up to the boundary for the Euler equations. {\it Nonlinearity}
{\bf 24}, 765--796.

\bibitem[Levermore \& Oliver(1997)]{LeOl}
{\sc Levermore, C.D.~\& Oliver, M.} 1997 Analyticity of solutions
for a generalized Euler equation. {\it J.~Diff. Equations} {\bf 133}, 321--339.

\bibitem[Lichtenstein(1925)]{Li25}
{\sc Lichtenstein, L.} 1925 \"Uber einige Hilfss\"atze der Potentialtheorie. I.
{\it Mathematische Zeitschrift\/} {\bf 23}, 72--88.

\bibitem[Lichtenstein(1927)]{Li27}
{\sc Lichtenstein, L.} 1927 \"Uber einige Existenzprobleme der Hydrodynamik.
{\it Mathematische Zeit\-schrift\/} {\bf 26}, 196--323.

\bibitem[Majda \& Bertozzi(2002)]{MB02}
{\sc Majda, A.J. \& Bertozzi A.L.} 2002 {\it Vorticity and incompressible flow}.
Cambridge Univ.~Press.

\bibitem[Matsubara(2008)]{Mat08}
{\sc Matsubara, T.} 2008 Resumming cosmological perturbations via the Lagrangian
picture: One-loop results in real space and in redshift space.
{\it Phys.~Rev.} D {\bf 77}, 063530.

\bibitem[Moutarde {\it et al.}(1991)]{Mou}
{\sc Moutarde, F., Alimi, J.M., Bouchet, F.R., Pellat, R.~\& Raman, A.} 1991
Precollapse scale invariance in gravitational instability. {\it Astrophys.~J.}
{\bf 382}, 377--381.

\bibitem[Pauls \& Matsumoto(2005)]{PaMa}
{\sc Pauls, W.~\& Matsumoto, T.} 2005 Lagrangian singularities of steady
two-dimensional flow. {\it Geophys.~Astrophys.~Fluid Dyn.} {\bf 99}, 61--75.

\bibitem[Pohle(1951)]{Phl}
{\sc Pohle, F.V.} 1951 The Lagrangian equations of hydrodynamics: solutions
which are analytic functions of the time. Thesis, New York University,
January 1951.

\bibitem[Rampf \& Buchert(2012)]{RB12}
{\sc Rampf, C.~\& Buchert, T.} 2012 Lagrangian perturbations and the
matter bispectrum I: fourth-order model for non-linear clustering
{\it J.~Cosmology and Astrop.~Phys.} JCAP06(2012)021.

\bibitem[Sahni \& Shandarin(1996)]{SS96}
{\sc Sahni, V.~\& Shandarin, S.} 1996 Accuracy of Lagrangian
approximations in voids. {\it Mon.~Not. R.~Astron.~Soc.} {\bf 282}, 641--645.

\bibitem[Serfati(1992)]{Se92}
{\sc Serfati, Ph.} 1992 \'Etude math\'ematique de flammes infiniment minces
en combustion. R\'esultats de structure et de r\'egularit\'e pour l'\'equation
d'Euler incompressible. Th\`ese de Doctorat de l'Universit\'e Paris 6.

\bibitem[Serfati(1995a)]{Se95a}
{\sc Serfati, Ph.} 1995a \'Equation d'Euler et holomorphies \`a faible r\'egularit\'e
spatiale. {\it C.~R.~Acad. Sci.~Paris.~S\'erie I\/} {\bf 320}, 175--180.

\bibitem[Serfati(1995b)]{Se95b}
{\sc Serfati, Ph.} 1995b Structures holomorphes \`a faible r\'egularit\'e spatiale
en m\'ecanique des fluides. {\it J.~Math.~Pures Appl.} {\bf 74}, 95--104.

\bibitem[Shnirelman(2012)]{Sh12}
{\sc Shnirelman, A.} 2012 On the analyticity of particle trajectories
in the ideal incompressible fluid. arXiv:1205.5837 [math.AP].

\bibitem[Stein(1970)]{Ste}
{\sc Stein, E.M.} 1970 {\it Singular integrals and differentiability
properties of functions}. Princeton Univ. Press.

\bibitem[Stoker(1957)]{Sto}
{\sc Stoker, J.J.} 1957 {\it Water waves: The mathematical theory
with applications}. Wiley-Interscience.

\bibitem[Stokes(1848)]{St48}
{\sc Stokes, G.G.} 1848 Notes on Hydrodynamics. IV. Demonstration
of a fundamental theorem. {\it Cambridge and Dublin Mathematical Journal}
{\bf III}, 209--219. Available at the Goettingen Archive.

\bibitem[Stokes(1883)]{St83}
{\sc Stokes, G.G.} 1883 Notes on Hydrodynamics. IV. Demonstration
of a fundamental theorem. In {\it Mathematical and physical papers}, vol.~II,
pp.~36--50. Cambridge Univ.~Press.

\bibitem[G.I.~Taylor(1938)]{Tay}
{\sc Taylor, G.I.} 1938 The spectrum of turbulence. {\it Proc.~R.~Soc.}
{\bf A164}, 476--490.

\bibitem[Taylor \& Green(1937)]{TG}
{\sc Taylor, G.I.~\& Green, A.E.} 1937 Mechanism of the production of small eddies
from large ones. {\it Proc.~R.~Soc.} {\bf A158}, 499--521.

\bibitem[Weber(1868)]{We68}
{\sc Weber, H.} 1868 Ueber eine Transformation der hydrodynamischen Gleichungen.
{\it Journal f\"ur die reine und angewandte Mathematik (Crelle)} {\bf 68},
286--292. Berlin.

\bibitem[Wolibner(1933)]{Wo33}
{\sc Wolibner, W.} 1933 Un th\'eor\`eme sur l'existence du mouvement plan
d'un fluide parfait, homog\`ene, incompressible, pendant un temps infiniment
long. {\it Mathematische Zeitschrift\/} {\bf 37}, 698--726.

\bibitem[Yakubovich \& Zenkovich(2001)]{YaZe}
{\sc Yakubovich, E.I.~\& Zenkovich D.A.} 2001 Matrix approach to Lagrangian
fluid dynamics. {\it J.~Fluid Mech.} {\bf 443}, 167--196.

\bibitem[Zheligovsky(2011)]{Zh11}
{\sc Zheligovsky, V.} 2011 A priori bounds for Gevrey--Sobolev norms of
space-periodic three-dimensional solutions to equations of hydrodynamic type.
{\it Advances in differential equations} {\bf 16}, 955--976.
arXiv:1001.4237 [math.AP]

\bibitem[Zygmund(2002)]{Zy}
{\sc Zygmund, A.} 2002 {\it Trigonometric series}. 3rd ed., Cambridge Univ.~Press.
\end{thebibliography}
\end{document}